\documentclass[amsmath,preprint,aps,fleqn]{revtex4}

\begin{document}

\title{Factorization of and Determinant Expressions for the Hypersums of Powers of Integers}
\author{Jerome Malenfant}
\affiliation{American Physical Society\\ Ridge, NY}
\date{\today}
\begin{abstract}
 We derive a compact determinant formula for calculating and factorizing the 
 hypersum polynomials $S^{(L)}_k(N) \equiv \sum_{n_1=1}^N \cdots 
 \sum_{n_{L+1}=1}^{n_{L}}(n_{L+1})^k$ expressed in the variable $N(N+L+1)$.  
\end{abstract}

\maketitle
\newtheorem{1}{Theorem}
\newtheorem{2}[1]{Theorem}
\newtheorem{3}[1]{Theorem}
\newtheorem{cor}{Corollary}
\newtheorem{lemma1}{Lemma}
\newtheorem{lemma2}[lemma1]{Lemma}
\newtheorem{I}{Proposition}

\section{Introduction}
In this article, we will consider the finite sums and hypersums of positive 
integers raised to a positive integer  $k$.  Let $S_{k}(N)$ denote the sum of 
such integers from 1 to $N$:
            \begin{equation} 
           S_{k}(N) =\sum^{N}_{n=1}n^k. 
           \end{equation}           
 The hypersums $S_k^{(L)}(N)$ are then defined recursively as
           \begin{eqnarray}
           S_k^{(0)}(N) = S_k(N),~~S^{(L)}_k(N) = \sum_{n=1}^N S_k^{(L-1)}(n)~~ {\rm for}~~L \geq 1.
           \end{eqnarray}
In particular,
           \begin{eqnarray}
           S^{(1)}_0(N) = \sum_{n=1}^N S^{(0)}_0(n) = \sum_{n=1}^N n 
		    = S^{(0)}_1(N),
           \end{eqnarray}
from which it follows that $S^{(L)}_0(N) = S^{(L-1)}_1(N)$ for all $L\geq 1$.  
In the following, we will mostly take $k>0$.  $S^{(L)}_k(N)$ is an 
$(L+k+1)$-order polynomial in $N$, and is given by the formula \cite{Knuth}
           \begin{eqnarray}
           S_k^{(L)} (N) = \sum_{q=0}^k S(k,q)~q! \left( \begin{array}{c}
		    N+L+1\\q+L+1 \end{array} \right)           
           \end{eqnarray}             
where $S(k,q)$ is  a Stirling number of the 2nd kind.   As an example, we 
quote from Knuth \cite{Knuth} a partial result of Faulhaber for $S^{(10)}_6(N)$ 
(which in Knuth's notation is $\Sigma ^{11}n^6)$: 
            \begin{eqnarray*}
            S_6^{(10)}(N) = \frac{5!}{17!} ~ \left\{ ~6~N^{17} + 561~N^{16} 
			+ \cdots + 1021675563656~N^5 + \cdots -96598656000~N ~
            \right\} .
            \end{eqnarray*}
                    
For $N=-1,-2, \ldots , -L-1$, the sum over $q$ on the right  side of 
eq.(4) is zero since the binomial coefficients are all zero;  the sum 
is also zero for $N=0$ if $k>0$, since in this case $S(k,0) =0$ and all 
$q>0$ binomial coefficients are zero.   The polynomial representing 
$S^{(L)}_k(N)$ for $k>0$ therefore has zeros at these values, and is 
expressible in the form
          \begin{eqnarray}
          S^{(L)}_k(N) = N(N+1)(N+2) \cdots (N+L+1) \times {\cal Q}^{(k-1)}(N)
          \end{eqnarray}
where ${\cal Q}^{(k-1)}(N)$ is a $(k-1)$-order polynomial in $N$.  

While the $S_k^{(L)}$ polynomial can be calculated using eq.(4), this  formula
does not readily lend itself to expression in the factored form of eq.(5).  In
addition, it is known that these polynomials simplify when expressed in the 
variable $N(N+L+1)$.  The aim of this article is to derive an alternative formula
for the $ S^{(L)}_k$ polynomial in this variable, one in which the zeros at
$N=0,-1, \ldots ,-L-1$ are factored out and in which the remaining $(k-1)$-order
polynomial is given as a determinant  of a relatively simple matrix, one that 
does not explicitly involve Stirling numbers (or any other ``complicated''
numbers).  We do this in Theorem 3 in Section IV; in particular, our result
for $S_6^{(10)}$ is  
     \begin{eqnarray*}
     ~\nonumber \\
     S_6^{(10)}(y) = \frac{6!}{17!} ~\frac{ \sqrt{ 4y+121}}{2} ~ y(y+10)(y+18) (y+24) (y+28)(y+30) \times \frac{3y^2+22y-220}{3},
     ~ \nonumber \\
     \end{eqnarray*}
where $y \equiv N(N+11)$.   To do this, we need to develop some machinery, so 
for now we will consider only the $L=0$ polynomials.  These are given by
Faulhaber's formula \cite{Wolfram2008}: 
            \begin{equation}
            S_k(N) = \frac{1}{k+1} \sum_{n=1}^{k+1}(-1)^{\delta_{nk}} \left( \begin{array}{c}
                  k+1\\
                  n
                  \end{array} \right)\; B_{k+1-n}N^n,
             \end{equation}
where $B_n$ is a Bernoulli number.   In the following section, we will derive an 
alternative expression for $S_k$ which seems to be more convenient in extending
the formalism to hypersums.

\section{Series Expansion}   

In the sum over $n$ for $S_{k+1}(N)$ we make the substitution $n \rightarrow m+1$, 
add and subtract  terms to make the sum on $m$ go from 1 to $N$,  and then expand 
$(m+1)^{k+1}$ binomially, with the result \cite{Wolfram2008, Guo1999}:
       \begin{equation}
       (k+1)S_k(N) = (N+1)^{k+1} -1 - \sum^{k-1}_{q=0} \left( \begin{array}{c} 
	   k+1\\q\end{array}\right)\;S_q(N).
       \end{equation} 
Repeating this procedure but now with the replacement $n \rightarrow m-1$, we 
get: 
    \begin{equation}
    (k+1)S_k(N) = N^{k+1} - \sum^{k-1}_{q=0} \left( \begin{array}{c}
           k+1\\
             q
    \end{array}\right)\;(-1)^{k-q}S_q(N).
    \end{equation} 
Adding (7) and (8) and dividing by $2(k+1)$, the odd-$q$ terms in the sum 
cancel if $k$ is an even integer, while for odd $k$ the even-$q$ terms 
cancel.  We therefore get separate recursion relations for the even and 
odd power sums:
 
            \begin{subequations} 
            \begin{eqnarray}
            S_{2p}(N)&=&  \frac{1}{2p+1}\left[ \frac {(N+1)^{2p+1}+ N^{2p+1}-1}{2} - \sum_{q=0}^{p-1} 
            \left( \begin{array}{c}
              2p+1\\
              2q
            \end{array}\right)
            S_{2q}(N) \right] ,\\
            S_{2p+1}(N)&=&\frac{1}{2p+2}\left[ \frac {(N+1)^{2p+2}+ N^{2p+2}-1}{2} -  \sum_{q=0}^{p-1} 
            \left( \begin{array}{c}
              2p+2\\
              2q+1
            \end{array}\right)
            S_{2q+1}(N)\right].
            \end{eqnarray}
            \end{subequations}
            
We now define $x$ as  $2N+ 1$ and the functions $f_p(x)$ and $g_p(x)$ as
            \begin{subequations}
            \begin{eqnarray}
            f_p(x)&=& \frac{1}{2^{2p+2}} \left[ (x+1)^{2p+1}+(x-1)^{2p+1}\right]
             -\frac{x}{2} ,\\
            g_p(x) &=& \frac{1}{2^{2p+3}} \left[ ( x+1)^{2p+2}+ (x-1)^{2p+2}\right] 
            - \frac{1}{2} -(p+1)~ \frac{x^2-1}{4}.
            \end{eqnarray}
            \end{subequations}
We will assume for now that $p > 0$.   Equations (9a) and (9b) then become
             \begin{subequations}
             \begin{eqnarray}
             S_{2p}(x) &=& \frac {1}{2p+1} \left[ f_p(x) - \sum_{q=1}^{p-1} 
             \left(
             \begin{array}{c}
               2p+1\\
               2q
             \end{array}\right)
             S_{2q}(x)
             \right],\\
             S_{2p+1}(x)&=&\frac{1}{2p+2} \left[ g_p(x)- \sum_{q=1}^{p-1}
             \left( \begin{array}{c}
               2p+2\\ 
               2q+1
             \end{array}\right)
             S_{2q+1}(x) \right].
             \end{eqnarray}
             \end{subequations}
(Note that we have absorbed the $q=0$ terms in the sums in (9a) and (9b) 
into $f_p$ and $g_p$, respectively.)  

$f_p(x)$ is an odd function of $x$, while $g_p(x)$ is even.    Further, 
$f_p(x)$ has zeros at $x=0, \pm 1$, and $g_p(x)$ has double zeros at 
$x=\pm 1$.  $S_{2p}(x)$ and $S_{2p+1}(x)$ are therefore odd and even 
functions of $x$, respectively, with zeros at these points. These sums 
can therefore be expanded in these functions, which we do as
            \begin{subequations}
            \begin{eqnarray}
            S_{2p} &=& (2p)!  \sum_{q=1}^{p} \frac{C_{p-q}^{(1)} }{(2q+1)!}~ f_{q},\\
            S_{2p+1}&=& (2p+1)! \sum_{q=1}^{p} \frac{C_{p-q}^{(2)}}{(2q+2)!}~ g_{q}.
            \end{eqnarray}
            \end{subequations}    
Equations (11) fix the values $C_0^{(1)}= C_0^{(2)}=1$.  We now substitute (12a)
into (11a) and (12b) into (11b) to get 
            \begin{subequations}
            \begin{eqnarray}
            \sum_{q=0}^p \frac{C_{p-q}^{(1)}}{(2q+1)!} = \sum_{q=0}^p \frac{ C_{p-q}^{(2)}}{(2q+1)!} =0 ~~~~(p>0).  \label{a}
            \end{eqnarray}
These two sets of coefficients satisfy the same recursion relation and 
have the same initializing value; they are therefore equal, and we will 
denote the $p$th coefficient simply as $C_p$.  These coefficients satisfy the additional recursion relations:
     \begin{eqnarray}
     &&\sum_{q=0}^p \frac{2^{2q}C_{p-q}}{(2q+1)!} = \frac{1}{(2p)!} ,   \label{b} \\
     &&\sum_{q=0}^p \frac{2^{2q+1}C_{p-q}}{(2q+2)!}  =  \frac{1}{(2p+1)!} ,  \label{c} \\
     &&\sum_{q=0}^p \frac{2^{2p-2q}C_{p-q}}{(2q)!} = C_p , \label{d} \\
     &&\sum_{q=0}^{p}  \frac{ 2^{2p-2q} C_{p-q}}{(2q+1)!} = \frac{E_{2p}}{(2p)!}, \label{e} \\
     &&\sum_{q=0}^{p}  \frac{E_{2q}C_{p-q}}{(2q)!} = 2^{2p} C_p,  \label{f} \\
     &&\sum_{q=0}^p \frac{C_{p-q}}{(2q)!} =  \frac{C_p}{2^{1-2p}-1}, \label{g} \\
     &&\sum_{q=0}^p \frac{1}{(2q+1)!} \frac{C_{p-q}}{2^{2p-2q}-2} = -\frac{ \delta_{0p}}{2}
             -\frac{1}{2(2p)!}, \label{h} \\
     &&\sum_{q=0}^p \frac{1}{(2q)!} \frac{C_{p-q}}{2^{2p-2q}-2}=  \frac{C_p}{2^{2p}-2}
            -\frac{1}{2(2p-1)!} ,  \label{i} \\
     &&\sum_{q=0}^p C_q C_{p-q} = \frac{2p-1}{1-2^{1-2p}}~ C_p ,  \label{j} 
     \end{eqnarray}       
and in general, for $N$ an integer $> 1$,                     
     \begin{eqnarray}
     &&\sum_{q=0}^p \frac{C_{p-q}}{(2q+1)!} ~N^{2q+1} = \frac{2}{(2p)!} \left[ (N-1)^{2p} + (N-3)^{2p} + \cdots + \left\{ \begin{array}{l} 
            1\\
            2^{2p} \end{array} \right.  \right]  \left\{ \begin{array}{l} 
            N ~{\rm even}\nonumber \\
            N~{\rm odd}~(p>0) \end{array} \right\}, \nonumber \\
            ~  \label{k} \\
            &&\sum_{q=0}^p \frac{C_{p-q}}{(2q)!}~ N^{2q}
             =  \left\{ \begin{array}{c}
            1\\
             (2^{1-2p}-1)^{-1} \end{array} \right\}  C_p \label{l} \\
             && ~~~~~~~~~~~~~~~~~~~~+  \frac{2}{(2p-1)!} \left[  (N-1)^{2p-1}
            + (N-3)^{2p-1} + \cdots +  \left\{ \begin{array}{l} 
            1\\
            2^{2p-1} \end{array} \right.  \right]  \left\{ \begin{array}{l} 
            N ~{\rm even}\\
            N~{\rm odd} \end{array} \right\} . \nonumber 
      \end{eqnarray}                        
      \end{subequations}
In (\ref{e}) and (\ref{f}), $E_k$ is the $k$th Euler number.   For the derivation of relations (\ref{b}-\ref{l}), 
see Appendix A.

Consider now the infinite sum $  \sum_{p=0}^{\infty} C_p x^{2p}$, which we assume to be convergent in 
some region.  From (13a) we have
       \begin{eqnarray}
       \sum_{p=0}^{\infty} C_p x^{2p}
       &=& 1-  \sum_{p=1}^{\infty}\sum_{q=0}^{p-1} \frac{ C_{q}}{(2p-2q+1)!} ~ x^{2p}  \nonumber \\
       &=& 1-  \sum_{q=0}^{\infty} C_q  \sum_{p=q+1}^{\infty} \frac{ x^{2p}}{(2p-2q+1)!} \nonumber \\
        &=& 1-  \sum_{q=0}^{\infty} C_q  \sum_{n=1}^{\infty} \frac{ x^{2q+2n}}{(2n+1)!}\nonumber  \\
       &=& 1-  \sum_{q=0}^{\infty} C_q x^{2q} \left(\frac{ \sinh x}{x} -1 \right) ,
                      \end{eqnarray} 
from which we get
      \begin{equation}
       \sum_{p=0}^{\infty} C_p x^{2p}= \frac{ x}{\sinh x} 
       \end{equation} 
As a consequence, the $C_p$ coefficients are related to the Bernoulli
numbers as 
     \begin{equation}
     C_p =  \frac {2- 2^{2p}}{(2p)!}~ B_{2p},  \label{31}
     \end{equation}
 and correspond to integer sequences A036280 (numerators) and A036281
(denominators) in the Encyclopedia of Integer Sequences \cite{Sloane2007}.  
The Bernoulli numbers are related to the Riemann zeta function 
$\zeta(s)$ for $s$ an even, positive integer (\cite{Edwards}, p. 12):
            \begin{equation}
             B_{2p} = (-1)^{p+1} \frac {2(2p)!}{(2 \pi)^{2p}}~ \zeta(2p).
             \end{equation}
Combining this with (\ref{31}),
            \begin{equation}
            C_{p} = (-1)^{p} \frac {2}{ \pi^{2p}} (1-2^{1-2p}) \zeta(2p).
            \end{equation}
But since             
            \begin{equation}
             \left(1-\frac{2}{2^{x}} \right)  \zeta(x)= \left(1-\frac{2}{2^{x}} \right)  \sum _{n=1}^{\infty} \frac{1}{n^x}
             = 1 - \frac{1}{2^x} + \frac{1}{3^x} - \frac{1}{4^x} + \frac{1}{5^x} - \cdots
             \end{equation}
we have              
     \begin{equation}          
     C_p = (-1)^p \frac{2}{\pi ^{2p}} \left[ 1 - \frac{1}{2^{2p}} + \frac{1}{3^{2p}} - \frac{1}{4^{2p}} + \frac{1}{5^{2p}} - \cdots \right]   = (-1)^p \frac{2}{\pi ^{2p}}~ \eta (2p),
     \end{equation}  
where $\eta (x)$ is the Dirichlet eta function.

Recursion relation (\ref{a}) corresponds to eq.(3.2) in Van Malderen 
\cite{Van}, and these $C_p$ coefficients are, up to a factor of $(-1)^p$, 
the same as Van Malderen's $D_p$ coefficients in that article; (see also the
article by Chen \cite{Chen}).  Relations (\ref{g}) and (\ref{l}), for $N=3$,
correspond to eqs. (2.2) and (2.3), respectively, in Van Malderen.

\begin{1}
For $k \geq 0$, $S_k(x)$ is given by the expression
           \begin{eqnarray*}
            S_{k}(x) =  \frac{k!}{2^{k+1}} \sum_{q=0}^{ \lfloor k/2 \rfloor } C_{q} ~
            \frac{  x^{ k+1-2q} -1}{(k+1-2q)!} .      
            \end{eqnarray*}			
\end{1}
$Proof.$  The sums in eqs.(12) for $S_{2p}$ and $S_{2p+1}$ can be extended down 
to $q=0$ since $f_0=g_0=0$.  The resulting sums are however not valid for 
the $p=0$ sums $S_0$ and $S_1$.  They can be made so by adding and subtracting 
terms inside the summation which sum to zero for $p>0$ as a consequence of
recursion relation (\ref{a}).  The ``corrected'' sums are:
            \begin{subequations}
            \begin{eqnarray}
            S_{2p}(x) &=& (2p)!  \sum_{q=0}^{p} \frac{C_{p-q} }{(2q+1)!} \left\{ \frac{ (x+1)^{2q+1}
             +(x-1)^{2q+1}} {2^{2q+2} } - \frac{1}{2} \right\} ,\\
            S_{2p+1}(x)&=& (2p+1)! \sum_{q=0}^{p} \frac{ C_{p-q}}{(2q+2)!}  \left\{ \frac{ (x+1)^{2q+2} 
            +(x-1)^{2q+2}} {2^{2q+3} } - \frac{1}{2} \right\} .
            \end{eqnarray}
            \end{subequations}              
Expanding the functions inside the brackets  in powers of $x$, we have:
            \begin{eqnarray}
             \frac{ (x+1)^{2q+1} +(x-1)^{2q+1}} {2^{2q+2} } - \frac{1}{2} &=&  \frac{1}{2^{2q+1}} \sum_{n=0}^q
            \left( \begin{array}{c}
            2q+1\\
            2n+1 \end{array} \right) x^{2n+1} - \frac{1}{2} \nonumber \\
            &=& \frac{1}{2^{2q+1}} \sum_{n=0}^q
            \left( \begin{array}{c}
            2q+1\\
            2n+1 \end{array} \right)( x^{2n+1} -1) ,
           \end{eqnarray}
and
            \begin{eqnarray}
             \frac{ (x+1)^{2q+2} +(x-1)^{2q+2}} {2^{2q+3} } - \frac{1}{2} 
             = \frac{1}{2^{2q+2}} \sum_{n=0}^q
            \left( \begin{array}{c}
            2q+2\\
            2n+2 \end{array} \right)( x^{2n+2} -1) .      
           \end{eqnarray}
Then,
           \begin{eqnarray}
            S_{2p}(x) &=& (2p)!  \sum_{q=0}^{p} \frac{C_{p-q} }{(2q+1)!}~ \frac{1}{2^{2q+1}} \sum_{n=0}^q
            \left( \begin{array}{c}
            2q+1\\
            2n+1 \end{array} \right)( x^{2n+1} -1) \nonumber \\
            &=&  (2p)! \sum_{n=0}^p \frac{ x^{2n+1} -1 } {(2n+1)!} \sum_{q=n}^{p} 
            \frac{2^{-2q-1} C_{p-q} }{(2q-2n)!}  .             
            \end{eqnarray}
The sum over $q$ is evaluated using (\ref{d}):
            \begin{eqnarray}
            \sum_{q=n}^p \frac{ 2^{-2q-1} C_{p-q} }{(2q-2n)!} = 
            \sum_{l=0}^{p-n} \frac{ 2^{-2n-1-2l} C_{p-n-l} }{(2l)!}  = \frac{  C_{p-n} }{2^{2p+1}},
            \end{eqnarray}
and so we have
             \begin{eqnarray}
             S_{2p}(x) = \frac{(2p)!}{2^{2p+1}} \sum_{n=0}^{p} C_{p-n} ~ \frac{  x^{2n+1} -1}{(2n+1)!}    
              =    \frac{(2p)!}{2^{2p+1}} \sum_{q=0}^{p} C_{q} ~ \frac{  x^{2p+1-2q} -1}{(2p+1-2q)!} . 
             \end{eqnarray}
Similarly, for $S_{2p+1}$:                      
             \begin{eqnarray}
              S_{2p+1}(x) &=&  (2p+1)! \sum_{q=0}^{p} \frac{C_{p-q} }{(2q+2)!}~ \frac{1}{2^{2q+2}}
              \sum_{n=0}^{q}
              \left( \begin{array}{c}
              2q+2\\
              2n+2 \end{array} \right) (x^{2n+2}-1)\nonumber \\
              &=& \frac{ (2p+1)!}{2}   \sum_{n=0}^p \frac{x^{2n+2}-1}{(2n+2)!} \sum_{q=n}^{p} 
              \frac{ 2^{-2q-1}C_{p-q} }{(2q-2n)!} .
              \end{eqnarray}
The sum over $q$ is the same as in (24) above, and so we get
            \begin{eqnarray}
             S_{2p+1}(x)= \frac{(2p+1)!}{2^{2p+2}} \sum_{n=0}^p C_{p-n}  \frac{x^{2n+2}-1 }{(2n+2)!}
              = \frac{(2p+1)!}{2^{2p+2}} \sum_{q=0}^p C_q ~\frac{x^{2p+2-2q}-1}{(2p+2-2q)!}.
             \end{eqnarray}
\\
QED

It is interesting to compare the expression in Theorem 1 with Faulhaber's formula. 
Converting to Bernoulli numbers and expressed in the variable $N$, it becomes.
      \begin{equation}
      S_k(N) = \frac{1}{2^k(k+1)} \sum_{n=1}^{k+1} \left( \begin{array}{c}
                k+1\\
                n       
		    \end{array} \right) (1- 2^{k-n}) B_{k+1-n} [(2N+1)^n-1],
      \end{equation}

The power-sum  polynomials take a simpler form in the variable $y=N(N+1)$ compared
to $x$.  (The fact that the polynomials for the odd-power sums simplified when
expressed in thid variable was known by Faulhaber;  see the discussion by Knuth 
\cite{Knuth}.)   In this variable, the polynomials are:             
             \begin{eqnarray}
             S_0 &=& ~ \frac{  \sqrt{4y+1}}{2}-\frac{1}{2},
             ~~~~~~~~~~~~~~~~~~~~~~~~~~~~~~~~ S_1 = \frac {y}{2}, \nonumber  \\
             S_2 &=& \frac { \sqrt{4y+ 1} ~y}{6}, 
            ~~~~~~~~~~~~~~~~~~~~~~~~~~~~~~~~~~~~S_3 = \frac {y^2}{4},\nonumber  \\
             S_4 &=& \frac { \sqrt{4y+1}~y(3y-1) }{30} , 
             ~~~~~~~~~~~~~~~~~~~~~~~~~S_5 = \frac {y^2(2y-1)}{12} ,\\
             S_6 &=& \frac { \sqrt{4y+1} ~y(3y^2-3y+1)}{52} ,
             ~~~~~~~~~~~~~~~~~S_7 = \frac {y^2(3y^2-4y+2)}{24} ,
             \nonumber \\
             S_8 &=& \frac { \sqrt{4y+1}~y(5y^3-10y^2+9y-3)}{90},
             ~~~~~~~S_9 =\frac {y^2(y-1)(2y^2-3y+3)}{20}, \nonumber \\
             &\vdots&~~~~~~~~~~~~~~~~~~~~~~~~~~~~~~~~~~~~~~~~~~~~~~~~~~~~~~~ \vdots
             \nonumber
             \end{eqnarray}          

\begin{cor}  The power sum of odd integers, 
            \begin{eqnarray*}
            {\bar S}_k(N) \equiv 1 + 3^k + 5^k + \cdots + N^k, ~~N = {\rm odd~integer} .
            \end{eqnarray*}   
is given by the expression,
            \begin{eqnarray*}
            {\bar S}_{k} (N) = \frac{k!}{2} \sum_{q=0}^{ \lfloor k/2 \rfloor } C_{q} ~ \frac{  (N+1)^{k+1-2q}}{(k+1-2q)!} . 
            \end{eqnarray*}
\end{cor}			
$Proof:$ This follows directly from relations (\ref{k}) and (\ref{l}).
     
Therefore, expressed in the variable 
            \begin{eqnarray}
            {\bar y} = \frac{N(N+2)}{4},
             \end{eqnarray}
the ${ \bar S}_k$ are, up to an overall constant and a factor of  $2^k$, given by the same
polynomials as in (30):
            \begin{eqnarray}
            {\bar S}_k ({\bar y}) = 2^k \left[ S_k({\bar y}) -  S_k \left(- \frac{1}{4} \right) \right] .
            \end{eqnarray}

\section{Matrix formalism}
Recursion relation (\ref{a}) can be written in triangular-matrix form as
            \begin{eqnarray}
             \left( \begin{array}{ccccc}
            1& ~&~&~&~ \\
            \frac{1}{3!} & 1 & ~& ~~\mbox{\Huge 0} &~\\
            \frac{1}{5!} & \frac{1}{3!} & 1&~&~  \\
            \vdots & ~ & ~& ~& \ddots
            \end{array} \right) \left( \begin{array}{c}
            C_0\\
            C_1\\
            C_2\\
            \vdots 
            \end{array} \right) = \left( \begin{array}{c}
            1\\
            0\\
            0\\
            \vdots 
            \end{array} \right) .
            \end{eqnarray}
Lower-triangular matrices which have 1's along their main diagonal are unit 
lower-triangular matrices.  Matrices which are constant along all diagonals are 
Toeplitz matrices.   The matrix above is thus a unit lower-triangular Toeplitz
(LTT) matrix.   

If $A$ is an $n \times n$ lower-triangular matrix, $(n \leq \infty)$, then, 
for $k<n$, the $k-truncation$ of $A$ is the $k \times k$ lower-triangular 
matrix obtained from $A$ by removing all rows and columns greater than $k$.  
It is straightforward to show that the product of two truncated 
lower-triangular matrices is equal to the truncation of the product and, 
consequently, the inverse of a truncated lower-triangular matrix is the 
truncation of the inverse matrix.

We define the infinite-dimensional lower shift matrix  $J$ as
           \begin{eqnarray}
            J  \equiv \left( \begin{array}{ccccc}
           ~0 ~& ~&~&~&~ \\
           ~1 ~& ~0~ & ~&  \mbox{\Huge 0} &~\\
           ~0 ~& ~1~ &~ 0 &~&~ \\
           \vdots & ~ & ~& ~& \ddots
           \end{array} \right);  ~~~~( J^p)_{ij} = \delta_{p,i-j};
            ~~~~{J^pJ^q = J^{p+q}}.
           \end{eqnarray}  
and, for future reference, the unit column vectors ${\bf I}_p$ 
           \begin{eqnarray}
           {\bf I}_1 \equiv
           \left( \begin{array}{c}
           1\\
           0\\
           0\\
           \vdots 
           \end{array} \right) ;       ~~~~  {\bf I}_{p} \equiv  {\it J}^{p-1} {\bf I}_1.
          \end{eqnarray}
Any infinite-dimensional LTT matrix $A$ can be expanded out in powers of 
$J$, with $J^0$ being the  identity matrix:
            \begin{eqnarray}
            A = \left( \begin{array}{ccccc}
            a_0 & ~&~&~&~ \\
            a_1 & a_0 & ~& ~~\mbox{\Huge 0} &~\\
            a_2 & a_1 & a_0 &~&~  \\
            \vdots & ~ & ~& ~& \ddots
            \end{array} \right) = \sum_{q=0}^{\infty} a_q  J^q. 
            \end{eqnarray}          
LTT matrices commute with one another since ${\it J^p}$ commutes with 
${\it J^q}$.  The determinant of a finite-dimensional lower-triangular 
matrix is the product of its diagonal elements, so all finite unit lower-triangular 
matrices have determinant 1.

We will be concerned in the following with determinants that have the general structure:
             \begin{eqnarray}
             \left| \begin{array}{cccccc}
              a_0 & ~ & ~ & ~& ~& c_0 \\
              a_1& a_0  & ~ & \mbox{\Huge 0} ~~ &~  & c_1  \\  
              \vdots & ~& \ddots &~& ~ & \vdots \\
              a_{k-2} &~&~&~& a_0  & c_{k-2} \\
              a_{k-1} & a_{k-2}  & \cdots &~&  a_1 & c_{k-1}  \end{array} \right| 
              \end{eqnarray}
To describe them, we introduce some notation:  The determinant above is of 
a $k \times k$ matrix consisting of the sum of a ``almost LTT'' matrix, 
one whose $(k,k)$ element is zero, (the 'base' matrix), and a matrix whose 
only non-zero elements are along the last column, (the 'tower' matrix).  For 
$0 \leq q \leq k-1$, we define the ``unit tower matrix'' $ K_{q}$ to be a 
$k\times k$ matrix whose elements are all zero except for the $(k-q,k)$
element, whose value is 1:
            \begin{eqnarray}
            K_{q} \equiv  \left( \left. \begin{array}{ccccc}
            ~ & ~  & ~ &~& ~    \\
            ~ & ~ & ~&~&~\\
            ~ & ~ & \mbox{\Huge 0}_{(k,k-1)}   &~& ~  \\          
            ~&~&~&~&~ \\
            ~ &  ~ & ~& ~ &~  \end{array} \right|
            \begin{array}{c} 
            ~0~\\
             \vdots\\
             ~1~ \\
             \vdots\\
            ~0~ \end{array} \right)  \begin{array}{c}
            ~\\
            ~\\
            ~ \\
             \left. \begin{array}{c} ~ \\  ~ \end{array} \right\} \\
             ~   \end{array}
             \begin{array}{c}
             ~\\
             ~\\
             ~ \\
             q\\
             ~   \end{array}. 
            \end{eqnarray}              
(The dimension $k$ of these matrices will usually be clear from the context.) 
Matrices such as $A$ in (36) and others defined below are infinite-dimensional,
whereas in taking determinants the matrices are finite-dimensional.  To keep
the notation simple, when taking determinants we will use the same symbols for 
the truncated, finite-dimensional matrices as the infinite-dimensional ones, 
but will indicate with a subscript  the dimensionality of the matrices inside 
the brackets.  In this notation, the determinant in (37) is     
     \begin{eqnarray}
     \det \left\{ A - a_0K_0   + \sum_{q=0}^{k-1} c_{k-1-q}  K_q  \right\}_{(k)} .    
     \end{eqnarray}  
We can also use this notation when $A$ is a non-Toeplitz, lower-triangular 
matrix with constant diagonal element $a_0$. 

Some useful identities are given in the following lemma:
\begin{lemma1}  
Let $A$ be a lower-triangular matrix which is constant 
($=a_0$) along the main diagonal, $T$ be an arbitrary tower matrix,  $L$ be a unit 
LTT matrix, and $D$ be a diagonal matrix:  $D=  {\rm diag}~\{ d_1,d_2, \ldots \}$.  
Then     
    \begin{eqnarray*}
    {\rm (I)}~~~~~~~~~~~~~~~~~~~ \det \{  A - a_0K_0   + x T \}_{(k)} &=& x 
	~  \det \{  A - a_0K_0   + T  \}_{(k)} ,\\          
    ~~~~{\rm (II)} ~~~~  \det \left\{  A - a_0K_0   + \sum_{q=0}^{k-1} c_{k-1-q}  K_q  \right\}_{(k)} &=&
      \sum_{q=0}^{k-1} c_{k-1-q}~\det \left\{ A - a_0K_0 + K_q \right\}_{(k)},\\  
   {\rm (III)}\! ~~~~~~~~~~~~~~~~~~ \det \{  A - a_0K_0   + L T \}_{(k)} &=& 
   \det \{ L^{-1} A - a_0K_0   +  T \}_{(k)} \\                   
   {\rm (IV)} ~~~~~~~~~~~~~~~~~  \det \{  A - a_0K_0 + DT \}_{(k)} &=& 
   d_k \det \{ D^{-1} A D- a_0K_0   +  T \}_{(k)}                             
             \end{eqnarray*}
\end{lemma1}
$Proof$: The first identity is a statement in this notation of the well-know 
property that a common factor of any column in a matrix can be factored 
out of its determinant.  The 2nd identity corresponds to an expansion 
by minors of the determinant on the left along the $k$th column.  The 
3rd and 4th identities follow from the identity $ (\det X)( \det Y)
= \det (XY)$, and the easily demonstrated relations ( $k \times k$ 
truncation assumed):  $L^{-1}K_0 =K_0,~ D^{-1}K_0D = K_0$, and $TD 
= d_kT$.

In particular, for tower matrices $T_1$ and $T_2$, we have from (II),
            \begin{eqnarray}
            \det \left\{ A-a_0K_0 + T_1 + T_2 \right\}_{(k)} =   \det \left\{ A- a_0K_0 + T_1 \right\}_{(k)} 
           +  \det \left\{ A- a_0K_0  + T_2 \right\} _{(k)} .
           \end{eqnarray} 
                 
We define the matrices $P$ and $Q$ as 
           \begin{subequations}         
           \begin{eqnarray}
           P &\equiv&  
          \left( \begin{array}{ccccc}
           ~1 ~& ~&~&~&~ \\
           ~ \frac{1}{3!} ~& ~1~ & ~&  \mbox{\Huge 0} &~\\
           ~ \frac{1}{5!}~& ~ \frac{1}{3!}~ &~ 1 &~&~ \\
           \vdots & ~ & ~& ~& \ddots
           \end{array} \right) =  \sum_{q=0}^{\infty} \frac{J^q}{(2q+1)!};\\
           Q &\equiv&  
          \left( \begin{array}{ccccc}
           ~1 ~& ~&~&~&~ \\
           ~ \frac{1}{2!} ~& ~1~ & ~&  \mbox{\Huge 0} &~\\
           ~ \frac{1}{4!}~& ~ \frac{1}{2!}~ &~ 1 &~&~ \\
           \vdots & ~ & ~& ~& \ddots
           \end{array} \right) =  \sum_{q=0}^{\infty} \frac{J^q}{(2q)!}.
           \end{eqnarray}   
           \end{subequations}                     
Then the $C_p$ coefficients are 
          \begin{eqnarray}
          C_p =  \det \left\{ P - K_0 + K_{p} \right\}_{(k)}
          \end{eqnarray}
for any $k >p$.  This is demonstrated by noting that $C_0=1$ by this formula and that
           \begin{eqnarray}
           \sum_{q=0}^p \frac{C_{p-q}}{(2q+1)!} &=& \sum_{q=0}^p \frac{1 }{(2q+1)!} \det \left\{ P - K_0 +
           K_{p-q} \right\}_{(k)} \nonumber \\
           &=&\det \left\{ P - K_0 +  \sum_{q=0}^p \frac{J^q}{(2q+1)!}   K_{p}    \right\}_{(k)}
           \end{eqnarray}   
is zero since the last column and the $(p+1)$ column in the determinant are equal.   It follows from
(42) that
            \begin{eqnarray}
            P^{-1} =  \sum_{q=0}^{\infty} C_q J^q 
            \end{eqnarray}
Written out, and expanding by minors, the RHS of eq. (42) is the $p \times p$ determinant
             \begin{eqnarray}
             C_{p} =   (-1)^{p} \left| \begin{array}{ccccc}
             \frac{1}{3!}  & 1 & 0 &  \cdots & 0 \\
            \frac{1}{5!} & \frac{1}{3!} & 1 &  \cdots & 0\\  
             \vdots & ~&~&~&\vdots \\
             \frac{1}{(2p-1)!} & \frac{1}{(2p-3)!} & \frac{1}{(2p-5)!}  & \cdots & 1\\
            \frac{1}{(2p+1)!} & \frac{1}{(2p-1)!} & \frac{1}{(2p-3)!}  & \cdots & 
             \frac{1}{3!} \end{array} \right| 
            \end{eqnarray}   
which is the form, up to the factor $(-1)^p$, that Van Malderen gives for his $D_p$ coefficients \cite{Van}.
 
 We now apply this formalism to the power sums $S_k(x)$ for $k \geq 2$:  
          \begin{eqnarray}
           S_{2p}(x) &=& \frac{(2p)!}{2^{2p+1}} \sum_{n=0}^{p} C_{p-n} ~ \frac{  x^{2n+1} }{(2n+1)!} ; ~~
           S_{2p+1}(x)= \frac{(2p+1)!}{2^{2p+2}} \sum_{n=0}^p C_{p-n}  \frac{x^{2n+2}-1 }{(2n+2)!}.          
           \end{eqnarray}
We now define the diagonal matrices $D_m$ as
             \begin{eqnarray}
            D_{m} \equiv
            \left( \begin{array}{ccccc}
            m!& ~ & ~& ~ & ~ \\
            ~ & (m+2)! & ~ &~ \mbox{\Huge 0} & ~ \\
            ~  & ~ &(m+ 4)! & ~ & ~ \\
            ~ & \mbox{\Huge 0}~~& ~ & \ddots& ~ 
            \end{array} \right) ,
            \end{eqnarray}
where $m$ takes on nonnegative integer values.  These matrices satisfy 
the identities
             \begin{eqnarray}
             D_m ^{ \pm 1}J^p {\bf I}_1& =& [(m+2p)!]^{ \pm 1}~ J^p {\bf I}_{1},
             \end{eqnarray}
as well as the corresponding identities for the truncated matrices, with ${\bf I}_1$ replaced by an 
appropriate unit tower matrix.
 
We have 
     \begin{subequations}
     \begin{eqnarray}
            S_{2p}(x)  &=& \frac{(2p)!}{2^{2p+1}} \sum_{n=0}^{p}\frac{  x^{2n+1} }{(2n+1)!} ~\det \left\{ P - K_0 + J^{n} 
            K_{p} \right\}_{(p+1)}\nonumber \\
           &=& \frac{(2p)!}{2^{2p+1}}~ \det \left\{ P - K_0 + \sum_{n=0}^{p}\frac{  x^{2n+1} }{(2n+1)!} ~J^{n} K_{p} \right\}_{(p+1)}
           \nonumber \\
           &=& \frac{(2p)!~x}{2^{2p+1}}~ \det \left\{ P - K_0 + D_1^{-1}\sum_{n=0}^{p}  x^{2n} J^{n} K_{p} \right\}_{(p+1)} \\ 
           &=& \frac{(2p)! ~x}{2^{2p+1}}~ \det \left\{ P - K_0 + D_1^{-1} ~  \frac{1}{I-x^2J} ~ K_{p} \right\}_{(p+1)} .\nonumber
      \end{eqnarray}       
And, in a like fashion,       
          \begin{eqnarray}
          S_{2p+1}(x)  &=& \frac{(2p+1)!}{2^{2p+2}}~ \sum_{n=0}^{p} \frac{x^{2n+2} -1}{(2n+2)!}~  \det
          \left\{ P - K_0 + J^{n} K_{p} \right\}_{(p+1)}\nonumber \\ 
          &=& \frac{(2p+1)!(x^2-1)}{2^{2p+2}}~  \det \left\{ P - K_0 + D_2^{-1} ~ \frac{1}{(I-J)(I-x^2J)} ~ K_{p} \right\}_{(p+1)}.
         \end{eqnarray}
         \end{subequations}
Note that       
         \begin{eqnarray}
         (I-J) D_1PD_1^{-1} {\bf I}_1 =  (I-J)^2 D_2PD_2^{-1} {\bf I}_1 = {\bf I}_1
          \end{eqnarray}   
I.e., the first column of $(I-J) D_1PD_1^{-1} $ and of 
$(I-J)^2 D_2PD_2^{-1}$ consists of 0's except for the $(1,1)$ element, 
which is 1.  As a consequence,       
         \begin{subequations}                       
         \begin{eqnarray}
         \det \left\{ P - K_0 + D_1^{-1} ~  \frac{1}{I-J} ~ K_{p} \right\}_{(p+1)} &=& \det \left\{ (I-J)D_1PD_1^{-1} - K_0 
         + K_{p} \right\}_{(p+1)} =    0 , \nonumber \\
         ~\\
         \det \left\{ P - K_0 + D_2^{-1} ~  \frac{1}{(I-J)^2} ~ K_{p} \right\}_{(p+1)}  &=& \det \left\{ (I-J)^2D_2PD_2^{-1} - K_0
         +K_p \right\} _{(p+1)} =0. \nonumber \\  
         ~
          \end{eqnarray}
          \end{subequations}
Using these identities, we can thus subtract off the tower matrices 
$D_1^{-1}(I-J)^{-1} K_p$ and $D_2^{-1}(I-J)^{-2} K_p$ from the tower 
matrices in (49a) and (49b), respectively, without changing the values 
of the determinants, to get, 
       \begin{subequations}
       \begin{eqnarray}
      S_{2p}(x)  &=& \frac{(2p)! ~x(x^2-1)}{2^{2p+1}}~ \det \left\{ P - K_0 + D_1^{-1} ~  \frac{J}{(I-J)(I-x^2J)} 
      ~ K_{p} \right\}_{(p+1)} ;   \\
       S_{2p+1}(x)  &=& \frac{(2p+1)! ~(x^2-1)^2}{2^{2p+2}}~ \det \left\{ P - K_0 + D_2^{-1} ~  \frac{J}{(I-J)^2(I-x^2J)} 
      ~ K_{p} \right\}_{(p+1)} .   
       \end{eqnarray} 
       \end{subequations}
The matrix $J$, acting on $K_p$ inside the brackets, lowers it to 
$K_{p-1}$.  The $(1,p+1)$ element of the matrix in each determinant is 
therefore zero, and we can expand each by minors along the first row, 
reducing them to $p \times p$ determinants.  Changing to the variable 
$y$ and using identities (III) and (IV) in Lemma 1, 
we get:     
            \begin{subequations}
            \begin{eqnarray}
            S_{2p}(x)  &=& \frac{y~\sqrt{4y+1}}{2^{2p-1}(2p+1)}~ \det \left\{ [ ~I-(4y+1)J~](I-J)D_3PD_3^{-1} - K_0 
            +  K_{p-1} \right\}_{(p)}; \\
            S_{2p+1}(x)  &=& \frac{y^2}{2^{2p-2}(2p+2)}~ \det \left\{ [~I-(4y+1)J~](I-J)^2D_4PD_4^{-1} - K_0 
            +  K_{p-1} \right\}_{(p)}   .    
            \end{eqnarray} 
            \end{subequations}
Theorem 3 in the next section generalizes these expressions to the $L>0$ 
cases.     
					       
\section{HYPERSUMS}
We  generalize relation (44) and define the coefficients 
$C_p^{(L)},~L \geq 0,$ by
        \begin{eqnarray}
       P^{-(L+1)} = \sum_{p=0}^{\infty} C_p^{(L)} J^p .
         \end{eqnarray}
We then have
         \begin{eqnarray}   
         C_p^{(0)} = C_p,~~~C_0^{(L)} =1,~~~C_p^{(L)} = \sum_{q=0}^p C_q^{(L-1)}C_{p-q}, ~~L \geq 1,
         \end{eqnarray}  
 and 
         \begin{eqnarray}
         C_{p} ^{(L)} = \det \left\{ P^{L+1} - K_0 +  K_{p} \right\}_{(k)} .
         \end{eqnarray}
If we multiply this determinant on the left by $ \det (P^{R-L-1})(=1)$, 
where $R$ is an arbitrary integer, and take the matrix $P^{R-L-1}$ inside, 
(56) becomes, with the replacement $p \rightarrow k-1-q$,
            \begin{eqnarray}
            C_{k-1-q} ^{(L)}  =   \det \left\{ P^{R} - K_0 + P^{R-L-1} J^q K_{k-1} \right\}_{(k)} .      
           \end{eqnarray}   
This second form is more useful in the proof of Theorem 3 below, where 
we will be adding terms containing $C^{(L)}_p$'s with different $L$ 
and $p$ values ; we can sum these up into a single determinant if the 
base matrices are the same.

We state below some recursion relations that these coefficients satisfy:
            \begin{subequations}       
            \begin{eqnarray}
            \sum_{q=0}^p \frac{2^{2q+1} C^{(L+2)}_{p-q}}{(2q+2)!} 
            &=& \sum_{q=0}^p \frac{ C^{(L+1)}_{p-q}}{(2q+1)!} =  C_p^{(L)} ;\label{x}\\
            \sum_{q=0}^p \frac{2^{2q}C^{(L+2)}_{p-q}}{(2q+1)!} &=& \sum_{q=0}^p
             \frac{C^{(L+1)}_{p-q}} {(2q)!}  = \frac{L+1-2p}{L+1} ~ C_p^{(L)}. \label{y}
             \end{eqnarray}
             \end{subequations}             
(The first equality in each line also holds for  $L=-1$.)   Using the third
equation in (55), the two equalities in (\ref{x}) follow from (\ref{c}) and 
(\ref{a}), respectively, and the first equality in (\ref{y}) follows from 
(\ref{b}).   To prove the second equality, it follows from (54) that 
the $C_p^{(L)}$'s are the coefficients in the 
expansion of $(x/\sinh x)^{L+1}$: 
         \begin{eqnarray}
           \left(  \frac {x}{\sinh x}\right)^{L+1} =  \sum _{p=0}^{\infty} C^{(L)}_p x^{2p}.
            \end{eqnarray}
Differentiating both sides             
            \begin{eqnarray}
            \frac {(L+1)x^{L}}{\sinh^{L+1} x} - \frac {(L+1)x^{L+1} \cosh x}{\sinh^{L+2} x} 
             =   \sum _{p=0}^{\infty} 2p~C^{(L)}_p x^{2p-1}.
            \end{eqnarray}
Rearranging terms and expanding out the $\cosh x$ factor, we get
            \begin{eqnarray}
            \sum _{p=0}^{\infty} \frac{ L+1-2p}{L+1} ~ C^{(L)}_p x^{2p} = 
            \frac {x^{L+2} \cosh(x)}{\sinh^{L+2}x} = \sum_{k=0}^{\infty} \sum_{q=0}^{\infty} C_k^{(L+1)}~
             \frac{x^{2k+2q}}{(2q)!}.
            \end{eqnarray}     
The proof then follows by equating powers of $x$. 

We define the variables $x_L$ and $y_L$:
           \begin{eqnarray}
            x_L= 2N+L+1; ~~ y_L=N(N+L+1);~~x_L^2= 4y_L+(L+1)^2.
            \end{eqnarray} 
(To keep the notation simple, we will in general suppress the subscript 
$L$ on these variables: $x_L \rightarrow x,~ y_L \rightarrow y$, with 
the $L$ understood.)  The generalization of Theorem 1 to hypersums is:  
 \begin{2} 
For $L \geq 0$ and $k>0$, the hypersum polynomials  are given by:  
              
              \begin{eqnarray*}   
             S_{2p}^{(2M)}(x)  &=&  \frac{(2p)!}{2^{2p+1+2M}} ~ \sum_{q=0}^{p+M} \left\{  C^{(2M)}_{p+M-q}~
             \frac{ x^{2q+1}}{(2q+1)!} -  \frac{1}{(2q+1)!} \sum_{s=1}^{M}  \frac{C_{p+M+1-q-s}^{(2M+2-2s)}}{(2s-1)!} ~
             \frac{x(x+2s-3)!!}{(x-2s+1)!!} ~ \right\} ; \\
                          ~\\
             S_{2p}^{(2M+1)}(x)  &=&  \frac{(2p)!}{2^{2p+2+2M}}~ \sum_{q=0}^{p+M} \left\{
              C^{(2M+1)}_{p+M-q}~ \frac{ x^{2q+2} }{(2q+2)!} -  \frac{1}{(2q+1)!} \sum_{s=1}^{M}
             \frac{ C_{p+M+1-q-s}^{(2M+2-2s)}}{(2s)!} ~\frac{x(x+2s-2)!!}{(x-2s)!!} ~ \right\} ;\\                      
             ~\\
              S_{2p+1}^{(2M-1)} (x) &=& \frac{(2p+1)!}{2^{2p+1+2M}} ~ \sum_{q=0}^{p+M}
              \left\{ C_{p+M-q}^{(2M-1)}~ \frac{x^{2q+1}}{(2q+1)!} - \frac{1}{(2q)!}  
              \sum_{s=1}^{M} \frac{C^{(2M-2s)}_{p+M+1-q-s}}{(2s-1)!} ~\frac{(x+2s-2)!!}{(x-2s)!!} \right\};\\
              ~\\             
             S_{2p+1}^{(2M)} (x) &=& \frac{(2p+1)!}{2^{2p+2+2M}} ~\sum_{q=0}^{p+M} 
             \left\{ C_{p+M-q}^{(2M)} ~\frac{x^{2q+2}-1}{(2q+2)!} - \frac{1}{(2q)!} \sum_{s=1}^{M} 
               \frac{C^{(2M-2s)}_{p+M+1-q-s}}{(2s)!} ~\frac{(x+2s-1)!!}{(x-2s-1)!!} \right\} ; \\      
              \end{eqnarray*}
\end{2}
where we use the convention that $C^{(L)}_p=0$ for $p<0$.  See Appendix B for the proof.

The sums over $q$ in the second terms in these expressions can be 
evaluated using some of the relations in (58), but the unsummed 
expressions above are more convenient in the proof of Theorem 3 
below.  Before getting to that theorem, we first need to prove 
a technical lemma:
\begin{lemma2}
The matrix $P$ satisfies the identities,                  
               \begin{eqnarray*}
               {\rm (I)}\!~~~~~~~ D_1 P^{2k+1}D_1^{-1}J^k~ {\bf I}_{1} &=&  \frac{J^k}{ (I-J) (I-9J) \cdots 
               [I-(2k+1)^2J] }~ {\bf I}_{1}, \\
               {\rm (II)} ~~~~~~ D_2 P^{2k+2}D_2^{-1}J^k~{\bf I}_{1} &=&  \frac{J^k}{ (I-4J) (I-16J) \cdots 
               [I-(2k+2)^2J ] }~ {\bf I}_{1}, \\   
               {\rm (III)} ~~~ D_1QP^{2k+1}D_1^{-1} J^k~ {\bf I}_{1} &=& 
               \frac{J^k}{ (I-4J) (I-16J) \cdots  [ I-(2k+2)^2J ]}~ {\bf I}_{1}, \\              
               {\rm (IV)} ~~~ D_2Q P^{2k+2}D_2^{-1}J^k~{\bf I}_{1} &=& 
                \frac{J^k}{ (I-J) (I-9J) \cdots [ I-(2k+3)^2J ] }~  {\bf I}_{1}.\\    
               \end{eqnarray*}
\end{lemma2}
$Proof.$   (We remind the reader that the $D_m$ matrices are defined by equation (47).)
Each equation above is a statement that the first column of the matrix 
on the left is equal to the first column of the matrix on the right.   The 
matrices themselves are not equal since the ones on the right are LTT and 
those on the left are not.  

Let the matrices $X_k$ and $Y_k$ be defined as:
             \begin{subequations}
            \begin{eqnarray}
             X_k &\equiv&  \frac{J^{k}}{ (I-J) (I-9J) \cdots [ I-(2k+1)^2J ] } ;\\
             Y_k &\equiv& \frac{J^{k}}{ (I-4J) (I-16J) \cdots [ I-(2k+2)^2J ] }.    
             \end{eqnarray}
             \end{subequations}
These matrices can be expanded, respectively, in powers of odd and even integers as
              \begin{subequations}
              \begin{eqnarray}
             X_k&=& \sum_{q=0}^k \frac {\alpha_q^{(k)}}{I-(2q+1)^2J }
              =\sum_{n=0}^{\infty} J^{n} \sum_{q=0}^k  \alpha_q^{(k)} (2q+1)^{2n} ,\\
             Y_k &=& \sum_{q=0}^k \frac {\beta_q^{(k)}}{I-(2q+2)^2J }
             =  \sum_{n=0}^{\infty} J^{n} \sum_{q=0}^k  \beta_q^{(k)} (2q+2)^{2n} .
             \end{eqnarray}
            \end{subequations}
The proofs of the first equalities in (64) follow by induction on $k$, and 
then the inverse matrices are expanded in positive powers of $J$.   It 
is easy to show that these sets of coefficients each satisfy $(k+1)$ 
independent equations,
             \begin{eqnarray}   
             \left. \begin{array}{l} 
              \sum^k_{q=0} \alpha_q^{(k)} (2q+1)^{2p}\\
              ~\\
              \sum^k_{q=0} \beta_q^{(k)} (2q+2)^{2p} \end{array} \right\} 
             =  \left\{ \begin{array}{l} 1, ~~{\rm for~} p=k;\\
             0, ~~{\rm for~}  0 \leq p < k; \end{array} \right.  \label{z}
             \end{eqnarray}
(project the far left and the far right sides of equations (64) onto ${\bf I}_1$:  
in powers of $J$, $X_k$ and $Y_k$ have the form  $J^k + a _1 J^{k+1}+ a_2 J^{k+2} 
+ \cdots$, so the resulting column vector's first $k$ elements are 0, with the 
$(k+1)$ elements equal to 1).  The expansions for $X_k$ and $Y_k$ are then:  
             \begin{subequations}           
             \begin{eqnarray}
             X_k &=& \sum_{n=0}^{\infty} J^{n+k} \sum_{q=0}^k 
             \alpha_q^{(k)} (2q+1)^{2n+2k} ;\\
              Y_k &= &\sum_{n=0}^{\infty} J^{n+k} \sum_{q=0}^k 
             \beta_q^{(k)} (2q+2)^{2n+2k}.
             \end{eqnarray}
             \end{subequations}
Relations (\ref{z}) are sufficient to determine the $(k+1)$ values of these 
coefficients.  The coefficients are in fact given by explicit formulas, but for 
the purpose of this proof, only these relations  are required.

We consider relation (I) first.  We will use the expansion for $P$ in 
powers of $J$ given in (41) and the $m=1$ identity in (48).  For $k=0$, 
the left-hand side of relation I reduces to
             \begin{eqnarray}
             D_1PD_1^{-1}~{\bf I}_1 = \sum_{q=0}^{\infty} J^q {\bf I}_1 = \frac{1}{I-J}~ {\bf I }_1,
             \end{eqnarray}
which equals $X_0 {\bf I}_1$ as required.   For $k>0$, an evaluation of the 
left-hand side of the relation results in
              \begin{eqnarray}
               D_1 P^{2k+1}D_1^{-1}J^k~ {\bf I}_1  = Z_k {\bf I}_1
               \end{eqnarray} 
where 
               \begin{eqnarray}
                Z_k &=& \sum_{n=0}^{\infty}\frac{ J^{k+n}}{(2k+1)!}
               \sum_{q_{2k}=0}^n \left( \begin{array}{c} 2n+2k+1\\2q_{2k} +2k \end{array} \right) \nonumber \\
               &~&~~~~ \times \sum_{q_{2k-1}=0}^{q_{2k}} \left( \begin{array}{c} 2q_{2k}+2k\\2q_{2k-1} +2k-1 
               \end{array} \right)  \cdots  \sum_{q_2=0}^{q_3} \left( \begin{array}{c} 2q_3 +3\\ 2q_2+2
                \end{array} \right) \sum_{q_{1}=0}^{q_{2}} \left( \begin{array}{c} 2q_{2} +2\\2q_{1} +1 
                \end{array} \right)   .
                \end{eqnarray}               
$Z_k$ $is$ an LTT matrix.  Relation I therefore holds if  and only if this matrix  
is equal to the matrix $X_k$ for all $k$.   We will assume then that, for some 
$k$, $Z_k$ is given by the sum
           \begin{eqnarray}
           Z_k  = \sum_{n=0}^{\infty} J^{n+k} \sum_{q=0}^k 
           \alpha_q^{(k)} (2q+1)^{2n+2k} .
           \end{eqnarray}
Then
              \begin{eqnarray}
             Z_{k+1} &=& \frac{(2k+1)!}{(2k+3)!} \sum_{n=0}^{\infty} J^{k+n+1} \sum_{p_{2}=0}^n 
             \left( \begin{array}{c}
              2n+2k+3\\2p_{2} +2k+2 \end{array} \right) \nonumber \\
              &~& ~~~~~~~~~~~~~~~~~~~~~~~~~ \times  \sum_{p_{1}=0}^{p_{2}} \left( \begin{array}{c} 
              2p_{2}+2k+2\\2p_{1} +2k+1
              \end{array} \right)  \sum_{q=0}^{k} \alpha_q^{(k)} (2q+1)^{2p_1+2k}.
              \end{eqnarray}
Using the binomial identities,            
         \begin{eqnarray*}
        \sum_{k=0}^{n} \left( \begin{array}{c} 2n+2\\2k+1 \end{array} \right) x^{2k+1} &=& \frac{(x+1)^{2n+2}
        - (x-1)^{2n+2}}{2} , \\
        \sum_{k=0}^{n} \left( \begin{array}{c} 2n+1\\2k \end{array} \right) x^{2k} &=& \frac{(x+1)^{2n+1}
        - (x-1)^{2n+1}}{2},
        \end{eqnarray*}  
and relations (\ref{z}) for the $\alpha$'s, this reduces to    
          \begin{eqnarray}
         Z_{k+1} &=& \sum_{n=0}^{\infty} J^{n+k+1} \sum_{q=0}^{k+1} (2q+1)^{2n+2k+2} \nonumber \\
          &~& ~~~~~~~~~~~~~~\times \frac{(2k+1)!}{4(2k+3)!} \left\{~  \frac{2q+1}{2q-1}~ \alpha_{q-1}^{(k)}
         -(2 + \delta_{0q})\alpha_q^{(k)}  + \frac{2q+1}{2q+3} ~\alpha_{q+1}^{(k)} 
         ~ \right\}; \\
          (\alpha_q^{(k)} &=&0 ~{\rm for ~} q<0 ~{\rm or~} >k).   \nonumber
          \end{eqnarray}
It is a straightforward calculation to show that the coefficients in the 
sum over $q$ satisfy relations (\ref{z}) for,  and thus must be equal to, 
$\alpha_q^{(k+1)}$.  As a consequence,  $Z_{k+1} = X_{k+1}$ and, by 
induction on $k$, for all $k'>k$.  

The proof of relation (II) follows by a similar argument, using the 
matrices $D_2$ and $Y_k$.
  
To prove (III),  we use relation (I) to write the left-hand side of (III) as 
            \begin{eqnarray}
            D_1QP^{2k+1}D_1^{-1}J^k~{\bf I}_1 = D_1QD_1^{-1} X_k ~{\bf I}_1 .
           \end{eqnarray}
Evaluating this using the expansion for $Q$ in (41b), we have:
            \begin{eqnarray}
            D_1QD_1^{-1} X_k ~{\bf I}_1
            &=&  D_1 \sum_{p=0}^{\infty} \frac{J^p}{(2p)!} D_1^{-1} \sum_{m=0}^{\infty} J^{m+k} 
           \sum_{q=0}^k  \alpha_q^{(k)} (2q+1)^{2m+2k}~ {\bf I}_1 \nonumber  \\
           &=& \sum_{n=0}^{\infty} J^{n+k} \sum_{q=0}^k \frac{ \alpha_q^{(k)}}{ 2(2q+1) } 
           \left\{ (2q+2)^{2n+2k+1}  +(2q)^{2n+2k+1} \right\} ~ {\bf I }_1\nonumber \\
           &=& \sum_{n=0}^{\infty} J^{n+k} \sum_{q=0}^k (2q+2)^{2n+2k} \left\{ \frac{q+1}{2q+1}
          ~\alpha_q^{(k)} +  \frac{ q+1}{2q+3} ~\alpha_{q+1}^{(k)} \right\}  {\bf I }_1  . 
          \end{eqnarray}      
Again, it can be shown that the coefficients in this equation satisfy 
(\ref{z}) for $\beta_q^{(k)}$, and we therefore have
           \begin{eqnarray}
            D_1QP^{2k+1}D_1^{-1}J^k~{\bf I}_1  = Y_k{\bf I }_1,
           \end{eqnarray}  
which is  relation (III).

The right-hand side of (IV) is 'almost' $X_{k+1}{\bf I}_1$, but it 
lacks one power of $J$.  Nevertheless, we can write the left-hand side, 
using (II), as
           \begin{eqnarray}
            D_2QP^{2k+2}D_2^{-1}J^k~{\bf I}_1 = D_2QD_2^{-1} Y_k ~{\bf I}_1,
           \end{eqnarray}
which evaluates as
            \begin{eqnarray}
            D_2QD_2^{-1} Y_k ~{\bf I}_1 
            &=& D_2 \sum_{p=0}^{\infty} \frac{J^p}{(2p)!} D_2^{-1} \sum_{m=0}^{\infty} J^{m+k} 
            \sum_{q=0}^k  \beta_q^{(k)} (2q+2)^{2m+2k} {\bf I}_1\nonumber \\
            &=& \sum_{n=0}^{\infty} J^{n+k} \sum_{q=0}^k
            \frac{ \beta_q^{(k)}}{ 2(2q+2)^2 } \left\{  (2q+3)^{2n+2k+2}+(2q+1)^{2n+2k+2} - 2 \right\} 
            {\bf I}_1  \nonumber \\
             &=& \sum_{n=0}^{\infty} J^{n+k} \sum_{q=0}^{k+1} (2q+1)^{2n+2k+2} 
            \left\{ \frac{ \beta_q^{(k)}} {2(2q+2)^2} +\frac{ \beta_{q-1}^{(k)}} {2(2q)^2} -\delta_{0q} d_k,
            \right\} {\bf I}_1    
            \end{eqnarray} 
where                                
           \begin{eqnarray*}
           d_k \equiv \sum_{q=0}^k \frac{ \beta_q^{(k)} }{(2q+2)^2} .
           \end{eqnarray*}
The coefficients in the sum over $q$ satisfy (\ref{z}) for $\alpha_q^{k+1}$ and 
therefore 
            \begin{eqnarray}
             D_2QP^{2k+2}D_2^{-1}J^k~{\bf I}_1  = \sum_{n=0}^{\infty} J^{n+k} \sum_{q=0}^{k+1} 
             \alpha_q^{(k+1)} (2q+1)^{2n+2k+2} {\bf I }_1,      
             \end{eqnarray} 
 which is relation (IV).  \\  
QED
    
We will use this Lemma in a slightly different form in proving Theorem 3 
below.  For ($k \times k$)-truncated matrices $J,~P$, etc, the identities 
in the Lemma also hold if the unit column vector ${\bf  I}_1$ is replaced 
by the $k \times k$ unit tower matrix $K_{k-1}$; both ${\bf I}_1$ and 
$K_{k-1}$ project out the first column of any matrix that they are multiplied 
by from the left.  We will refer to both forms (the ``column-vector'' and 
the ``tower-matrix'' forms) simply as ``Lemma 2''.

Our main result is:
\begin{3}
Expressed in the variable $y= N(N+L+1)$, the hypersum 
 polynomials $S_k^{(L)}(y)$ for $L \geq 0,~k>0$ are given by 
              \begin{eqnarray*}
              S^{(L)}_{k} (y) &=& \frac{k!}{(L+k+1)! } ~\frac{ ( \sqrt{~4y+(L+1)^2~})^{A}} {2^{k-1}(1+L~ {\rm mod}~2) }
              ~ \prod _{q=0}^{\lfloor L/2 \rfloor }  [ ~y+q(L+1-q)~ ] \\
              & \times& \det \left\{ [~I-(4y+(L+1)^2)J~] ~ \prod _{q=0}^ {\lfloor L/2 \rfloor } [~I-(L+1-2q)^2 J~ ] ~
               D_mP^{L+1}D_m^{-1}  - K_0  +  K_{n-1}  \right\}_{(n)}     
              \end{eqnarray*}              
where
            \begin{eqnarray*}
            A &=& L ~{\rm mod}~2 +(k+1)~{\rm mod}~2= 0,1,2;\\
            m &=& L+3- k~{\rm mod}~2;\\
            n &=& \left\lfloor \frac{k+1}{2} \right\rfloor = p,~p+1.\\
            \end{eqnarray*}
\end{3} 

$Proof.$   We consider first the $S^{(2M)}_{2p}$ polynomials.  Then 
$A=1$,  $m= 2M+3$, and $n=p$.  The sum over $q$ from 0 to $p+M$ in 
the expression for this polynomial in Thm II requires that the size 
of the matrix used to represent the coefficients be at least $(p+M+1)$.  
Using (57) with $R=2M+1$, we can express the $C^{(2M)} _{p+M-q}$ 
and the $C^{(2M+2-2s)} _{p+M+1-q-s}$ coefficients as determinants 
with the same base matrix:
              \begin{subequations}
              \begin{eqnarray}
              C^{(2M)} _{p+M-q} &=&  \det  \left\{ P^{2M+1} - K_0 +  J^{q} K_{p+M}\right\}_{(p+M+1)} ; \\ 
               C_{p+M+1-q-s}^{(2M+2-2s)} &=& \det  \left\{ P^{2M+1} - K_0 + P^{2s-2} J^{q+s-1} K_{p+M} \right\}_{(p+M+1)} .
              \end{eqnarray}  
              \end{subequations}      
Then we have from Thm 2 and identity (II) in Lemma 1 that
            \begin{eqnarray}
            S^{(2M)}_{2p}(x)  &=&  \frac{(2p)!}{2^{2p+1+2M}} ~ \sum_{q=0}^{p+M} \left\{  C^{(2M)}_{p+M-q}~
             \frac{ x^{2q+1}}{(2q+1)!} -  \frac{1}{(2q+1)!} \sum_{s=1}^{M}  \frac{C_{p+M+1-q-s}^{(2M+2-2s)}}{(2s-1)!} ~
             \frac{x(x+2s-3)!!}{(x-2s+1)!!} ~ \right\} \nonumber  \\           
             &=&  \frac{(2p)!}{2^{2p+1+2M}} \det  \left\{ P^{2M+1} - K_0 + \sum_{q=0}^{p+M} 
            \frac{ x^{2q+1}}{(2q+1)!}~J^{q} K_{p+M} \right.  \nonumber \\ 
            &-&  \left. \sum_{q=0}^{p+M} \frac{1}{(2q+1)!} \sum_{s=1}^{M} \frac{ P^{2s-2} J^{q+s-1}}{(2s-1)!} 
            ~ \frac{x(x+2s-3)!!}{(x-2s+1)!!} ~K_{p+M} ~ \right\} _{(p+M+1)}
             \end{eqnarray}
This expression will be reduced and simplified by making a series of 
replacements for terms and factors in the tower-matrix part which are 
valid inside the brackets of the determinant, (but not generally outside 
of them).  The first replacement is for the sum over $q$ in the 4th term:  
$\sum_q J^q/(2q+1)! \rightarrow P$.  Although $P$ is given by an infinite 
sum, this replacement is allowed since $P$ inside the brackets is 
truncated.  Other factorials in the denominators can be dealt with by 
making the additional replacements:
            \begin{eqnarray*}
            && \sum_{q=0}^{p+M} \frac{ x^{2q+1}}{(2q+1)!}~J^{q} K_{p+M} \rightarrow D_1^{-1} ~\frac{x}{I-x^2J}~K_{p+M}; \\
            && \frac{J^{s-1}}{(2s-1)!}~ K_{p+M} \rightarrow D_1^{-1}J^{s-1} K_{p+M}.
            \end{eqnarray*}
Eq.(80) then becomes                      
            \begin{eqnarray}
            S^{(2M)}_{2p}(x) &=& \frac{(2p)!}{2^{2p+1+2M}} ~\det  \left\{ P^{2M+1} - K_0 +D_1^{-1} \left( \frac{x}{I-x^2J}
              \right. \right. \nonumber \\
            &-&  \left. \left. \sum_{s=1}^{M}   D_1 P^{2s-1}D_1^{-1} J^{s-1}~ \frac{x (x+2s-3)!!}  {(x-2s+1)!!} ~\right) K_{p+M} 
            ~ \right\} _{(p+M+1)} .
            \end{eqnarray}
The ratio of double factorials is equal to $x$ for $s=1$ and is
            \begin{eqnarray}
            \frac{x (x+2s-3)!!}  {(x-2s+1)!!} = x(x^2-1)(x^2-9) \cdots (x^2- (2s-3)^2)
            \end{eqnarray}
for $s>1$.  We  can also make the replacement                                          
              \begin{eqnarray*}
              D_1P^{2s-1}D_1^{-1}J^{s-1} ~ \rightarrow ~\frac{J^{s-1}}{(I-J) (I-9J) \cdots  [I-(2s-1)^2J]}.
            \end{eqnarray*}             
This is allowed by the ``tower-matrix'' form of relation (I) in Lemma 2. 
Upon making these changes, we get
               \begin{eqnarray}
              S_{2p}^{(2M)}(x) &=& \frac{(2p)!}{2^{2p+1+2M}} ~ \det \left\{ P^{2M+1} 
              - K_0 +D_1^{-1} \left(  \frac{x} {I-x^2J} - \frac{x}{I-J}  \begin{array} {c} ~\\~ \end{array}  \right. \right.  \nonumber  \\
               && \left. \left.  \begin{array} {c} ~\\~ \end{array} -   \frac{x}{I-J} \sum_{s=2}^M~ \prod_{q=1}^{s-1} \frac{x^2-(2q-1)^2}
               {I-(2q+1)^2J} ~J^{s-1} \right) K_{p+M} \right\}_{p+M+1}
                \end{eqnarray}
The terms inside the  brackets can be summed as
             \begin{eqnarray}
              \frac{x} {I-x^2J} - \frac{x}{I-J}  -   \frac{x}{I-J} \sum_{s=2}^M~ \prod_{q=1}^{s-1} \frac{x^2-(2q-1)^2}
               {I-(2q+1)^2J} ~J^{s-1}
                = x \prod_{q=0}^{M-1} \frac{ x^2-(2q+1)^2}{I-(2q+1)^2J} ~  \frac{J^M}{I-x^2J}   \nonumber \\     
                \end{eqnarray}
by successively combining terms with the same power of $J$:
                \begin{eqnarray*}
                &&   \frac{x} {I-x^2J} - \frac{x}{I-J} = \frac{x(x^2-1)}{I-J} ~ \frac{J}{I-x^2J}, \\
                && \frac{x(x^2-1)}{I-J} ~\frac{ J}{I-x^2J} -  \frac{x(x^2-1) J}{(I-J)(I-9J)} 
                 = \frac{x(x^2-1)(x^2-9)}{(I-J)(I-9J)} ~ \frac{ J^2}{I-x^2J} , \\
                && {\rm etc.},
                \end{eqnarray*}
Inserting (84)  into eq. (83), factoring out the $x \prod (x^2-(2q+1)^2)$ (as per 
identity (I) in Lemma 1), and using $J^MK_{p+M} =K_p$ (inside the brackets), 
we have           
               \begin{eqnarray}
               S_{2p}^{(2M)}(x) &=&  \frac{(2p)!}{2^{2p+1+2M}}  ~x \prod _{q=0}^{M-1} [~ x^2-(2q+1)^2~ ]  \nonumber\\    
               & \times&  \det \left\{ P^{2M+1} - K_0 +D_1^{-1} \prod _{q=0}^{M-1} \frac{1}{I-(2q+1)^2J} 
               ~  \frac{1}{I-x^2J}~K_{p } \right\}_{(p+M+1)}.
               \end{eqnarray}

Now consider the value of the determinant in this equation at $x= 2M+1$.   
By relation (I) in Lemma 2,  this is equal to  
             \begin{eqnarray}
             \det \left\{ P^{2M+1} - K_0 +\frac{P^{2M+1}}{(2M+1)!} ~K_{p} \right\}_{(p+M+1)}. 
             \end{eqnarray}
This determinant however is zero, since it contains two columns, (the ($M+1$)-th 
and the last column), which are proportional to one another.  We can 
therefore make the replacement 
               \begin{eqnarray*}
                \frac{1}{I-x^2J} \rightarrow  \frac{1}{I-x^2J} - \frac{1}{I-(2M+1)^2 J} 
                = \frac{ x^2-(2M+1)^2 }{I-(2M+1)^2J} ~ \frac{J}{I-x^2J}          
               \end{eqnarray*}   
in the determinant in (85) without changing its value.  Factoring out 
$x^2-(2M+1)^2$, this equation then becomes:
             \begin{eqnarray}
               S_{2p}^{(2M)}(x) &=&  \frac{(2p)!}{2^{2p+1+2M}}  ~x \prod _{q=0}^{M} [~ x^2-(2q+1)^2~ ] \nonumber\\    
               & \times&  \det \left\{ P^{2M+1} - K_0 +D_1^{-1} \prod _{q=0}^{M} \frac{1}{I-(2q+1)^2J} 
               ~  \frac{1}{I-x^2J}~K_{p-1 } \right\}_{(p+M+1)}.
              \end{eqnarray}
                           
This expression accomplishes, for these values of $L$ and $k$, the goal 
set out in the Introduction, to derive a determinant formula for the 
hypersum polynomials which factors out the zeros at $N=0, \ldots,-L-1$.  
The size of the determinant can however be reduced; the first $(M+1)$ 
elements in the last column of the tower matrix  in (87) are zero and 
we can make the replacement in that equation of
                 \begin{eqnarray}
               D_1^{-1} \rightarrow \left( \begin{array} {c|c}
               ~~\mbox{\huge 0}_{(M+1)}~ & ~ \\ \hline
               ~ & ~D_{2M+3}^{-1}~ \end{array} \right),
               \end{eqnarray} 
Successively expanding the determinant by minors along the top row, (87) is 
thus reduced to the $p \times p$ determinant:  
               \begin{eqnarray}
               S_{2p}^{(2M)}(x) &=&  \frac{(2p)!}{2^{2p+1+2M}}  ~x \prod _{q=0}^{M} [~ x^2-(2q+1)^2~ ] \\    
               & \times&  \det \left\{ P^{2M+1} - K_0 +D_{2M+3}^{-1} \prod _{q=0}^{M} \frac{1}{I-(2q+1)^2J} 
               ~  \frac{1}{I-x^2J}~K_{p-1 } \right\}_{(p)} \nonumber.      
               \end{eqnarray}
Now applying identites (III) and (IV) in Lemma 1 and changing the 
variable to $y$ using $x^2= 4y +(L+1)^2$, eq. (89) takes the form as 
stated in the Theorem for $L=2M,~k=2p$.
                                  
Likewise, the expressions in Theorem 2 for the other hypersums can be written 
in determinant form as:
              \begin{eqnarray*}
              S^{(2M+1)}_{2p}(x) &=& \frac{(2p)!}{2^{2p+2+2M}} ~\det \left\{ P^{2M+2} -K_0 +D_2^{-1}
               \left( \frac{x^2}{I-x^2J} \right. \right.\\
              &-&  \left. \left. \sum_{s=1}^{M}   D_2 P^{2s}D_2^{-1} J^{s-1}~\frac{x(x+2s-2)!!}{(x-2s)!!} \right) K_{p+M} 
              ~ \right\} _{(p+M+1)} 
              \end{eqnarray*}
                \begin{eqnarray*}                                                                                                     
                 S^{(2M-1)}_{2p+1}(x) &=& \frac{(2p+1)!}{2^{2p+1+2M}} ~\det  \left\{ P^{2M} - K_0 +D_1^{-1}
             \left( \frac{x}{I-x^2J} \right. \right. \\
            &-&  \left. \left. \sum_{s=1}^{M}   D_1Q P^{2s-1}D_1^{-1} J^{s-1}~\frac{(x+2s-2)!!}{(x-2s)!!} \right) K_{p+M} 
            ~ \right\} _{(p+M+1)} 
            \end{eqnarray*}
             \begin{eqnarray*} 
             S^{(2M)}_{2p+1}(x) 
             &=& \frac{(2p+1)!}{2^{2p+2+2M}} ~\det  \left\{ P^{2M+1} - K_0 +D_2^{-1} \left( \frac{x^2-1}{(I-J)(I-x^2J)} \right. \right. \\
            &-&  \left. \left. \sum_{s=1}^{M}   D_2 QP^{2s-2}D_2^{-1} J^{s-1}~\frac{(x+2s-1)!!}{(x-2s-1)!!} \right) K_{p+M} 
            ~ \right\} _{(p+M+1)} 
            \end{eqnarray*}
Using relations (II), (III), and (IV), respectively, in Lemma 2, they can be 
brought into the forms stated in the Theorem by a similar process.  For 
$S_{2p+1}^{(2M)}$ and $S_{2p+1}^{(2M-1)}$, the determinants corresponding 
to that in (86) for $S_{2p}^{(2M)}$ are not equal to zero, and so the final 
forms for these polynomials contain the tower matrix $K_p$ rather than 
$K_{p-1}$.  \footnote {Except for the $L=0$, $k=2p+1>1$ case, which we've 
seen can be reduced to a $p \times p$ determinant, a result of the identity 
$(I-J)^2D_2PD_2^{-1} {\bf I}_1 = {\bf I}_1$.}\\
QED. 
                       
We have already given the expression for $S_6^{(10)}$ in the Introduction; 
some other examples are:    
            \begin{eqnarray*}
            S_6^{(3)}(y) &=& \frac{6!}{10!} ~y(y+3)(y+4) \times (y^2-2y-1) ; \\
            S_7^{(5)}(y) &=& \frac{7!}{13!} \sqrt{y+9}~ y(y+5)(y+8) \times \frac{7y^3+14y^2-238y+295}{7} ; \\
            S^{(8)}_{11}(x)  &=&  \frac{11!}{20!}  ~y(y+8)(y+14)(y+18)(y+ 20) \times  \frac{14y^5 -4011y^3+25868y^2
             +145896y - 1199616 }{14} \\
            S^{(14)}_{14}(y) &=&  \frac{14!}{29!}  ~\frac{ \sqrt{4y+225}}{2} ~y(y+14)(y+26)(y+36)(y+ 44)(y+ 50) (y+54)(y+56)\\    
               && ~~~~~~~~ \times  ( ~ y^6 -1750y^4 +29960y^3 + 376167 y^2 -11436860 y +62455917~ )
               \end{eqnarray*}
                
Now let $\Delta^{(L)}_{k}$ be the determinant in the expression for $S^{(L)}_k$:
              \begin{eqnarray*}
              \Delta ^{(L)}_{k} (y) \equiv \det \left\{ [~I-(4y+(L+1)^2)J~] ~ \prod _{q=0}^ {\lfloor L/2 \rfloor } [~I-(L+1-2q)^2 J~ ] ~
               D_mP^{L+1}D_m^{-1}  - K_0  +  K_{n-1}  \right\}_{(n)}    
                \end{eqnarray*} 
where $m$ and $n$ are as defined in Theorem 3.
              
\begin{I}
$\Delta^{(L)}_{k}$ has the series expansion
     \begin{eqnarray*}
     \Delta ^{(L)}_{k} (y) &=& 4^n  \sum_{s=1}^{n} y^{n-s} ~\frac{a_{s}}{4^s} 
     \end{eqnarray*}
where $a_s$ is the $s \times s$ determinant
              \begin{eqnarray*}
                a_s  &=&  \det \left\{ [I-(L+1)^2J]^{n+1-s} \prod _{q=0}^ {\lfloor L/2 \rfloor } [I-(L+1-2q)^2 J ] ~
                  D_{L+k+3-2s}P^{L+1}D_{L+k+3-2s}^{-1} 
                 - K_0  +  K_{s-1}  \right\}_{(s)} . \\
               \end{eqnarray*} 
\end{I}			   				    
$Proof:$ We re-express $\Delta^{(L)}_k$ as           
               \begin{eqnarray}
              \Delta ^{(L)}_{k} (y) &=& \det \left\{  \prod _{q=0}^ {\lfloor L/2 \rfloor } [~I-(L+1-2q)^2 J~ ] 
               D_mP^{L+1}D_m^{-1}  - K_0  +  \frac{1}{I-[4y+(L+1)^2]J  } ~K_{n-1}  \right\}_{(n)}. \nonumber \\  
                 \end{eqnarray}   
We have                  
              \begin{eqnarray}
               \frac{1}{I-[4y+(L+1)^2]J  } &=&  \sum_{q=0}^{\infty} [4y+(L+1)^2]^{q}J^q   
               = \sum_{q=0}^{\infty }  \sum_{s=0}^q \left( \begin{array}{c} q\\s \end{array}\right)  (4y)^s (L+1)^{2q-2s} J^q ,    
                 \end{eqnarray}    
although only the $q \leq n-1$ terms in the $q$-sum will contribute to the determinant.   We interchange the 
order of the sums and take the sum over $s$ outside of the determinant:
              \begin{eqnarray}
              \Delta ^{(L)}_{k} (y) &=&\sum_{s=0}^{n-1} (4y)^s  \det \left\{  \prod _{q=0}^ {\lfloor L/2 \rfloor }
               [~I-(L+1-2q)^2 J~ ] D_mP^{L+1}D_m^{-1}  - K_0  \right. \nonumber \\
               && \left. ~~~~~~~~~~~~~~~~~~~~~~ +  \sum_{q=s}^{n-1} \left( \begin{array}{c} q\\s \end{array}\right)
                (L+1)^{2q-2s}J^qK_{n-1} \right\}_{(n)}   .
                 \end{eqnarray}   
Using the negative binomial series,
               \begin{eqnarray}
               \sum_{r=0}^{\infty} \left( \begin{array}{c} r+s\\s \end{array}\right) z^r =
               \left( \frac{1}{1-z} \right)^{s+1},
               \end{eqnarray}       
the sum over $q$ becomes  
               \begin{eqnarray}
               \sum_{q=s}^{n-1} \left( \begin{array}{c} q\\s \end{array}\right) (L+1)^{2q-2s}J^qK_{n-1}
                &=& \sum_{q=s}^{\infty} \left( \begin{array}{c} q\\s \end{array}\right) (L+1)^{2q-2s}J^qK_{n-1} \nonumber \\
                 &=& \left( \frac{1}{I-(L+1)^2J} \right)^{s+1}K_{n-1-s}.
               \end{eqnarray}    
We now make the change of summation index  $s \rightarrow n-s$  and get    
              \begin{eqnarray}
              \Delta ^{(L)}_{k} (y) =   \sum_{s=1}^{n} (4y)^{n-s} a_{s} ,
              \end{eqnarray}
where
               \begin{eqnarray*}
                a_s  &=&  \det \left\{ [I-(L+1)^2J]^{n+1-s} \prod _{q=0}^ {\lfloor L/2 \rfloor } [I-(L+1-2q)^2 J ]
                D_{m}P^{L+1}D_{m}^{-1} - K_0  +  K_{s-1}  \right\}_{(n)} . 
               \end{eqnarray*}    
As was done in the proof of Theorem 3, the size of the determinant can 
be reduced from $n \times n$ to $ s \times s$ by the replacements 
$D_m^{ \pm 1} \rightarrow D_{m+2n-2s}^{ \pm 1} $.    We recall from Theorem 3
that  $m= L+3-k~{\rm mod}~2$ and $n = \lfloor (k+1)/2 \rfloor$; then 
$m+2n-2s$ is equal to $L+k+3-2s$ and the result follows.\\
 QED                                       
                              
\section{Conclusion}

We have derived a formula for the hypersum $S_{k}^{(L)}$ polynomial in a 
factored form which does not explicitly involve either Bernoulli or Stirling 
numbers.  The hard part of the calculation has been reduced to the 
multiplication of matrices constructed entirely from 'simple' numbers and to 
the calculation of a determinant of size $\lfloor (k+1)/2 \rfloor$.
 
 \begin{acknowledgments}
I wish to thank Yonko Millev for invaluable assistance.
\end{acknowledgments}

\appendix{}
\section{ Recursion relations} 

As noted above, recursion relation (\ref{a}) follows directly from inserting either 
(12a) into (11a) or (12b) into (11b).  The demonstration of most of the remaining 
relations is more transparent (or at least more compact) if these are written in matrix 
form.  Matrices  ${\it P}$, ${\it Q}$, and ${\it J}$ and the column vectors ${\bf I}_p$ 
were previously defined in section III;   we also define the diagonal matrix,
          \begin{eqnarray}
          {\it 2}  \equiv
            \left( \begin{array}{ccccc}
            2^0& ~ & ~& ~ & ~ \\
            ~ & 2^2 & ~ &~ \mbox{\Huge 0} & ~ \\
            ~  & ~ & 2^4 & ~ & ~ \\
            ~ & \mbox{\Huge 0}~~& ~ & 2^6 & ~ \\
            ~ &~ &  ~  &~ & \ddots 
            \end{array} \right)  ,
            \end{eqnarray}
and the column vectors,
          \begin{eqnarray}
           {\bf C} \equiv \sum_{q=0}^{\infty} C_q {\bf I}_{q+1} ,~~~
          {\bf {\bar C}} \equiv \sum_{q=0}^{\infty} {\bar C}_q {\bf I}_{q+1} ,~~~
          {\bar C}_q \equiv  \frac {E_{2q}}{(2q)!} ,
          \end{eqnarray}
The matrix elements of ${\it P^{-1}}$ are proportional to the Bernoulli numbers 
and those of ${\it Q^{-1}}$ are proportional to the Euler numbers:
                  \begin{eqnarray}
                  {\it P^{-1}} &=& \sum_{q=0}^{\infty}  C_q {\it J^q}  =  \sum _{q=0}^{\infty} 
                  \frac {2-2^{2q}}{(2q)!} ~B_{2q} {\it J^q} , \\                
                  {\it Q^{-1}} &=& \sum_{q=0}^{\infty} {\bar C}_q {\it J^q} =  \sum_{q=0}^{\infty} \frac 
                  {E_{2q}}{(2q)!} {\it J^q} .          
                  \end{eqnarray}
In this matrix notation, (\ref{a}) is
                 \begin{equation}
                 {\bf C} = {\it P}^{-1} {\bf I}_1 .
                  \end{equation}

Now consider the product $ PQ$:
                \begin{eqnarray}
                 PQ = \sum _{p,q =0}^{\infty} \frac{ J^{p+q}}{(2p+1)!(2q)!}
                &=& \sum_{n=0}^{\infty}  J^n~ \sum_{q=0}^{n} \frac{1}{(2n-2q+1)!(2q)!}\nonumber\\
                &=& \sum_{n=0}^{\infty} \frac{J^n}{(2n+1)!} \sum_{q=0}^{n}\left( \begin{array}{c}
                      2n+1\\
                      2q \end{array} \right).
               \end{eqnarray}  
The evaluation of the sum over $q$ is straightforward and equals $2^{2n}$, so we have
            \begin{equation}
             PQ  = \sum_{n=0}^{\infty} \frac{ 2^{2n} J^n}{(2n+1)!}.
            \end{equation}
In similar fashion, the products $P^2$ and $ Q^2$ can be evaluated:
             \begin{subequations}
             \begin{eqnarray}
              P^2  &=& \sum_{n=0}^{\infty} \frac{2^{2n+1} J^n}{(2n+2)!},\\
              Q^2 &=&  J^0 + \sum_{n=0}^{\infty} \frac{ 2^{2n+1} J^{n+1} }{(2n+2)!} = I + JP^2.
             \end{eqnarray}
             \end{subequations}   
On the other hand,              
            \begin{equation}
            {\it 2}J^n {\it 2} ^{-1} = 2^{2n} J^n, 
            \end{equation}
and therefore
            \begin{equation}
          {\it 2}P{\it 2}^{-1}  = \sum_{n=0}^{\infty} \frac{ {\it 2J^n2}^{-1}}{(2n+1)!} = PQ
            \end{equation}  
Also,
            \begin{eqnarray}
             JP^2  = \sum_{n=0}^{\infty} \frac{2^{2n+1} J^{n+1}}{(2n+2)!} 
            = \sum_{p=1}^{\infty} \frac{2^{2p-1} J^{p}}{(2p)!}
            = \frac{1}{2} \left\{ \sum_{p=0}^{\infty} \frac{2^{2p} J^{p}}{(2p)!} - I \right\},
            \end{eqnarray}
or
            \begin{equation}
           {\it 2Q2}^{-1} = {\it I} + 2 \cdot  JP^2, ~~(= 2 \cdot Q^2 -I ~{\rm from~(A8b)}).
            \end{equation}            
            
Now relation (\ref{b})
            \begin{equation}
            \sum_{q=0}^p \frac{2^{2q}C_{p-q}}{(2q+1)!} = \frac{1}{(2p)!} 
            \end{equation}
corresponds to the matrix equation
            \begin{eqnarray}
            {\it 2P2}^{-1} {\bf C} = Q {\bf I}_1 
            \end{eqnarray}
and this is derived by multiplying both sides of the equation ${\bf C} = P^{-1} {\bf I}_1$ 
on the left by $PQ(=QP)$ and using relation (A10).   The third relation, (\ref{c}), 
             \begin{eqnarray}
             \sum_{q=0}^p \frac{2^{2q+1}C_{p-q}}{(2q+2)!}   =  \frac{1}{(2p+1)!}
             \end{eqnarray}
is equivalent, under the substitutions $ {\bar p} =  p+1$ and $ {\bar q} = q+1$, (and then dropping the `bar' notation), to           
            \begin{eqnarray*}
            \sum_{q=0}^p \frac{2^{2q}C_{p-q}}{(2q)!} = C_p +  \frac{2}{(2p-1)!}.
            \end{eqnarray*}
In this form,  (\ref{c}) is (A12) above multiplied on the right by ${\bf C}$:  
            \begin{equation}
            {\it 2Q2}^{-1} {\bf C} = {\bf C}  + 2 \cdot JP {\bf I}_1.
            \end{equation}            
           
The next three recursion relations
            \begin{subequations}
            \begin{equation}
            (\ref{d}):~~  \sum_{q=0}^p \frac{2^{2p-2q}C_{p-q}}{(2q)!} = C_p ,
            \end{equation}
            \begin{equation}
            (\ref{e}):~~  \sum_{q=0}^{p}  \frac{ 2^{2p-2q} C_{p-q}}{(2q+1)!} = \frac{E_{2p}}{(2p)!},
            \end{equation}
            \begin{equation}
           (\ref{f}) :~~  \sum_{q=0}^{p}  \frac{E_{2q}C_{p-q}}{(2q)!} = 2^{2p} C_p,  
            \end{equation}   
            \end{subequations} 
all follow from relation (A10).  In matrix notation, (\ref{d}) is   
            \begin{eqnarray} 
            ({\it 2}^{-1}Q {\it 2}) {\bf C} = {\it 2}^{-1} {\bf C}, ~~ {\rm or} ~~  Q {\it 2} {\bf C} ={\bf C} .
            \end{eqnarray}
(\ref{d}) is proven by noting that the left-hand side of the equation on the right is
            \begin{equation}
            Q{\it 2} {\bf C} =( P^{-1}{\it 2}P {\it 2}^{-1}){\it 2} ( P^{-1} {\bf I}_1) =  P^{-1}{\it 2} {\bf I}_1 
             =  P^{-1} {\bf   I}_1 = {\bf  C},
             \end{equation}
(since  $ {\it 2} {\bf I}_1 =  {\bf I}_1 $). Now inverting and rearranging equation (A10), 
we have
             \begin{equation}
             {\it 2P^{-1}2^{-1}P} =  {\it Q^{-1}}.
             \end{equation}                               
The left-hand side of this equation is                 
             \begin{equation} 
            {\it 2P^{-1}2}^{-1}P = \sum_{p=0}^{\infty} J^p \sum_{q=0}^p 2^{2p-2q}
             \frac {C_{p-q}}{(2q+1)!},
             \end{equation}
while the right-hand side is given by (A4).  Equating powers of $J$,  we have
             \begin{equation}
             \sum_{q=0}^{p}  \frac{ 2^{2p-2q} C_{p-q}}{(2q+1)!} = {\bar C}_p ,
             \end{equation}
which is relation (\ref{e}).   Rearranging the factors in  (A20), it becomes
             \begin{equation}
             {\it 2}^{-1}Q^{-1}P^{-1}{\it 2} =  P^{-1} .
              \end{equation}
From (A3) and (A4), this is  relation (\ref{f}), apart from an overall factor of $2^{-2p}$:
              \begin{eqnarray}
               \sum _{q=0}^p 2^{-2p}{\bar C}_{q} C_{p-q} = C_p .
              \end{eqnarray}
              
To prove relation (\ref{g}), 
            \begin{equation}
            \sum_{q=0}^p \frac{C_{p-q}}{(2q)!} = \frac{C_p}{2^{1-2p}-1},
            \end{equation}    
consider the matrix
            \begin{eqnarray}        
            {\cal M } \equiv 2 \cdot  P{\it 2}^{-1}P^{-1}Q  -  Q.
             \end{eqnarray}
The right-hand side can be rewritten as below:
             \begin{eqnarray}        
            {\cal M } &=& 2 \cdot ({\it 2}^{-1}PQ) P^{-1}Q  -  Q ,~~ ({\rm  using} ~P{\it 2}^{-1} 
              = {\it 2}^{-1}QP) \nonumber \\
            &=&   2 \cdot {\it 2}^{-1}Q^2  -  Q  \nonumber \\
            &=& {\it 2}^{-1} ({\it 2}Q{\it 2}^{-1} + I) -Q ,~ ~({\rm using~} 2 \cdot Q^2 = {\it 2Q2}^{-1} + I )   
            \nonumber \\
            &=& ( Q +I ){\it 2}^{-1}  -Q .
             \end{eqnarray}  
Therefore,
            \begin{equation}
            {\cal M} {\bf I}_1 = (Q + I) {\it 2}^{-1}{\bf I}_1  -Q {\bf I}_1
              =  Q{\bf I}_1+  {\bf I}_1- Q{\bf I}_1 = {\bf I}_1;
              \end{equation}
 or 
            \begin{eqnarray}        
            {\bf I}_1 =  (2 \cdot  P{\it 2}^{-1}P^{-1}Q  -  Q) {\bf I}_1.
             \end{eqnarray} 
Multiplying both sides of this equation by ${\it P}^{-1}$, it becomes
            \begin{equation}
            {\bf C}  = (2 \cdot {\it 2}^{-1} - {\it I}) {\it Q} {\bf C}.
            \end{equation}
Solving for ${\it Q}{\bf C}$,  we get        
             \begin{eqnarray}
             {\it Q} {\bf C} = \frac{1}{ 2 \cdot {\it 2}^{-1} - {\it I} }~{\bf C }  ,     
              \end{eqnarray}
which is (\ref{g}).  We now write the inverse matrix in this equation as            
             \begin{eqnarray} 
             \frac{1} { 2 \cdot {\it 2}^{-1} - {\it I}} = - {\it I} - \frac{2}{ {\it 2} - 2 \cdot {\it I}}.
             \end{eqnarray} 
Equation (A31) then becomes                        
             \begin{eqnarray}
             \frac{ 1} {{\it 2} - 2 \cdot  {\it I}} {\bf C } = - \frac{1}{2} ( { {\bf C }  + \it Q} {\bf C} ) .   
             \end{eqnarray}
Recursion relations (\ref{h}),
            \begin{equation}
             \sum_{q=0}^p \frac{1}{(2q+1)!} \frac{C_{p-q}}{2^{2p-2q}-2}
              = - \frac{\delta_{0p}}{2} -\frac{1}{2(2p)!},
             \end{equation}
and (\ref{i}), 
            \begin{equation}
            \sum_{q=0}^p \frac{1}{(2q)!} \frac{C_{p-q}}{2^{2p-2q}-2}= \frac{C_p}{2^{2p}-2}
            -\frac{1}{2(2p-1)!}  ,
            \end{equation}
are obtained by multiplying this equation by ${\it P}$ and by ${\it Q}$.  Multiplying 
by ${\it P}$ we get          
              \begin{eqnarray}
               {\it P} \frac{1}{ {\it 2} - 2 \cdot {\it I}} {\bf C} =  -\frac{1}{2} ( {\it P} {\bf C}+  {\it PQ} {\bf C} ) 
               =  -\frac{1}{2} ({\it I + Q} ){\bf I}_1
              \end{eqnarray}
which is (\ref{h}).   If we multiply (A33) instead by ${\it Q}$, we have,
              \begin{eqnarray}
              {\it Q} ({\it 2} - 2 \cdot  {\it I})^{-1} {\bf C } 
              = -\frac{1}{2} ( {\it Q} {\bf C } + {\it Q}^2 {\bf C} ) .  
              \end{eqnarray}
 From (A8b), ${\it Q^2} {\bf C} = ({\it I + JP^2}) {\bf C} = {\bf C} + {\it JP } {\bf I}_1$, so 
              \begin{eqnarray}
              {\it Q}( {\it 2}- 2 \cdot {\it I})^{-1} {\bf C} 
              =  -\frac{1}{2} ( {\it Q} {\bf C} +  {\bf C}  + {\it JP} {\bf I}_1 ) 
               =  ({\it 2} - 2 \cdot {\it I})^{-1} {\bf C} - \frac{1}{2} {\it JP} {\bf I}_1,
              \end{eqnarray}
which is (\ref{i}).

By eq. (54), $P^{-2} {\bf I}_1 = {\bf C}^{(1)}$, where ${\bf C}^{(1)}$ is 
the column vector whose elements are the $C^{(1)}_p$ coefficients.  Therefore, 
multiplying Eq.(A29) by $P^{-2}$, we get                
                \begin{eqnarray}
                {\bf C}^{(1)} =  2 \cdot P^{-1} {\it 2}^{-1}P^{-1}Q {\bf I}_1  -Q {\bf C}^{(1)}
                \end{eqnarray}
or           
               \begin{eqnarray}     
               Q {\bf C}^{(1)} &=& 2 \cdot (P^{-1} {\it 2}^{-1}Q) P^{-1} {\bf I}_1  -{\bf C}^{(1)} \nonumber \\
               &=& 2 \cdot ( {\it 2}^{-1} P^{-2}) {\bf I}_1 - {\bf C}^{(1)}~~( {\rm using ~}             
               P^{-1} {\it 2}^{-1}Q = {\it 2}^{-1} P^{-1}~{\rm from~ (A20)}) \nonumber \\
                &=& (2 \cdot {\it 2}^{-1} -I ) {\bf C}^{(1)}
                \end{eqnarray}
which gives us
                \begin{eqnarray} 
                 \sum_{q=0}^p \frac{C_{p-q}^{(1)} }{(2q)!} = (2^{1-2p}-1) C_p^{(1)}
                 \end{eqnarray}
Comparing this relation to the second equality in (\ref{y}) for $L=0$, we get
                \begin{eqnarray}
                C^{(1)}_p = \frac{2p-1}{1-2^{1-2p}} ~C_p
                \end{eqnarray}
which is (\ref{j}).
 
To prove relations (\ref{k}) and (\ref{l}), consider the polynomial
            \begin{eqnarray}
            \sum_{q=0}^p \frac{C_{p-q}}{(2q+1)!} ~u^{2q+1}.
            \end{eqnarray}            
Using relation (\ref{31}), this is expressible in terms of Bernoulli polynomials,
            \begin{eqnarray}
            \sum_{q=0}^p \frac{C_{p-q}}{(2q+1)!}~ u^{2q+1} = \frac{2}{(2p+1)!} \left[ B_{2p+1}(u) 
            - 2^{2p}B_{2p+1}\left( \frac{u}{2} \right) \right].
            \end{eqnarray}  
As a consequence of the identity $B_k(u+1) =B_k(u) + ku^{k-1}$,  the RHS becomes 
             \begin{eqnarray}            
              \sum_{q=0}^p \frac{C_{p-q}}{(2q+1)!} ~u^{2q+1} &=&  \frac{2}{(2p)!} ~(u-1)^{2p} 
              +\sum_{q=0}^p \frac{C_{p-q}}{(2q+1)!} ~(u-2)^{2q+1} \nonumber \\
               &=& \frac{2}{(2p)!} \left\{( u-1)^{2p} + (u-3)^{2p}
            + \cdots + (u-2N+1)^{2p} \right\} \nonumber\\
            &~&~~~~~~~ + \sum_{q=0}^p \frac{C_{p-q}}{(2q+1)!} ~(u-2N)^{2q+1} ;
             \end{eqnarray}      
Differentiating this equation with respect to $u$, 
            \begin{eqnarray}
            \sum_{q=0}^p \frac{C_{p-q}}{(2q)!} ~u^{2q} &=&  \frac{2}{(2p-1)!} ~(u-1)^{2p-1} +
            \sum_{q=0}^p \frac{C_{p-q}}{(2q)!} ~(u-2)^{2q} \nonumber \\
             &=& \frac{2}{(2p-1)!} \left\{ (u-1)^{2p-1} + (u-3)^{2p-1}
            + \cdots + (u-2N+1)^{2p-1}\right \} \nonumber\\
            &~&~~~~~~~~~~~~+ \sum_{q=0}^p \frac{C_{p-q}}{(2q)!} ~(u-2N)^{2q} .
            \end{eqnarray}           
Setting $u=2N$ and $u=2N+1$ in these two sets of equations prove, respectively, relation (\ref{k}) and, using 
(\ref{g}), relation (\ref{l}).
                                                
\section{Proof of Theorem 2}    
The proof is by induction on $L$ and uses the expressions for the integer power sums 
from Theorem 1,
            \begin{subequations}            
            \begin{eqnarray}
             S_{2p}(N)&=& \frac{(2p)!}{2^{2p+1}} \sum_{q=0}^p C_{p-q}  \frac{(2N+1)^{2q+1}-1 }{(2q+1)!}
             \nonumber \\
               & =&  \frac{(2p)!}{2^{2p+1}} \sum_{q=0}^p C_{p-q}  \frac{(2N+1)^{2q+1} }{(2q+1)!} 
               ~~( {\rm for~ }p>0),\\           
             S_{2p+1}(N)&=& \frac{(2p+1)!}{2^{2p+2}} \sum_{q=0}^p C_{p-q} \frac{(2N+1)^{2q+2}-1 }{(2q+2)!},
             \end{eqnarray}
and the odd-integer power sums from Corollary 1,
            \begin{eqnarray}
            {\bar S}_{2p}(N) &=& \frac{(2p)!}{2} \sum_{q=0}^{p} C_{p-q} ~
            \frac{  (N+1)^{2q+1} }{(2q+1)!},  \\   
            {\bar S}_{2p+1}(N) &=& \frac{(2p+1)!}{2} \sum_{q=0}^{p} C_{p-q} ~
            \frac{  (N+1)^{2q+2}} {(2q+2)!} ,     
            \end{eqnarray}
            \end{subequations}  
where in the last two equations $N$ is an odd integer.   The expressions for 
$S^{(L)}_k$ in Theorem 2 reduce to the integer power sums above when $L=0$.  It 
is necessary then to show that these expressions satisfy 
                \begin{eqnarray}
                 \sum_{n=1}^N S_{k}^{(L)}(n)  = S_{k}^{(L+1)}(N).
                 \end{eqnarray}
Since $S_k^{(L)}(n) =0$ for $n=0,-1,-2, \ldots , -L-1$, we can make the replacement
                 \begin{eqnarray}
                \sum_{n=1}^N S_k^{(2M,2M+1)}(n)~~ \rightarrow ~~ \sum_{n=-M}^NS_k^{(2M,2M+1)}(n)
                \end{eqnarray}
without affecting the values of the sums.  Taking the sum from $n=-M$ to $n=N$ of the 
first term in the expression for $S^{(2M)}_{2p}$ in Theorem 2, we have
                 \begin{subequations}
                  \begin{eqnarray}
                \sum_{n=-M}^N  \sum_{q=0}^{ p+M} C^{(2M)}_{p+M-q}
                 \left.  \frac{ x^{2q+1} }{(2q+1)!} \right|_{x=2n+1+2M} &=& \sum_{q=0}^{ p+M} \frac {C^{(2M)}_{p+M-q}}{(2q+1)!}~
                 {\bar S}_{2q+1}(2N+1+2M)  \nonumber \\
                 &=& \frac{1}{2}  \left. \sum_{q=0}^{ p+M}  \sum_{s=0}^q 
                  C^{(2M)}_{p+M-q}   C_{q-s}~ \frac{ X^{2s+2} }{(2s+2)!} \right|_{X=2N+2+2M} \nonumber \\
                 &=& \frac{1}{2}  \left. \sum_{s=0}^{p+M}  \frac{X^{2s+2}}{(2s+2)!} 
                 ~ \sum_{q=s}^{p+M}C^{(2M)}_{p+M-q}C_{q-s} \right|_{X=2N+2+2M} \nonumber \\
                 &=& \frac{1}{2} \sum_{s=0}^{p+M}  C^{(2M+1)}_{p+M-s}~
                 \left. \frac{ X^{2s+2}}{(2s+2)!} \right|_{X=2N+2+2M}
                 \end{eqnarray}  
In the 2nd line we've used (B1d) for the sum of odd integers raised to an odd power.   
By similar calculations, we get for the other sums,               
               \begin{eqnarray}
                 \sum_{n=-M}^N \sum_{q=0}^{p+M} C^{(2M+1)}_{p+M-q}~\left. \frac{ x^{2q+2}}{(2q+2)!}\right|_{x=2n+2+2M}
                &=& \frac{1}{2}   \left.  \sum_{s=0}^{p+M+1} C^{(2M+2)}_{p+M+1-s}~\frac{X^{2s+1}-X}{(2s+1)!} \right|_{X=2N+3+2M}\\   
                 \sum_{n=-M+1}^N \sum_{q=0}^{p+M} C^{(2M-1)}_{p+M-q}~\left. \frac{ x^{2q+1}}{(2q+1)!} \right|_{x=2n+2M} 
                 &=& \frac{1}{2} \left.  \sum_{s=0}^{p+M}C^{(2M)}_{p+M-s}~  \frac{ X^{2s+2}-1}{(2s+2)!} \right|_{X=2N+1+2M}\\                                  
                 \sum_{n=-M}^N \sum_{q=0}^{p+M} C^{(2M)}_{p+M-q}~\left. \frac{ x^{2q+2}-1}{(2q+2)!} \right|_{x=2n+1+2M}  
                 &=& \frac{1}{2} \left.  \sum_{s=0}^{p+M+1} C^{(2M+1)}_{p+M+1-s}~ \frac{ X^{2s+1} -2^{2s}X}{(2s+1)!} 
                 \right|_{X=2N+2+2M}                                                             
                 \end{eqnarray}  
                 \end{subequations}

We now consider the 2nd terms in the expressions for $S_k^{(L)}$.   We can use the identities 
                \begin{eqnarray*}
                \left. \frac{ (x+2s-1)!!}  {(x-2s-1)!!}\right|_{x=2n+1+2M}
               &=& 2^{2s}(2s)!   \left( \begin{array}{c} n+M+s\\2s \end{array} \right),\\
               \left.  \frac{ (x+2s-2)!!}  {(x-2s)!!}\right|_{x=2n+2+2M}
               &=& 2^{2s-1}(2s-1)!   \left( \begin{array}{c} n+M+s\\2s -1\end{array} \right),
               \end{eqnarray*}                               
               \begin{eqnarray*}
               \sum_{j=s}^N \left( \begin{array}{c} j\\s \end{array} \right) = \left( \begin{array}{c} N+1\\s+1
               \end{array} \right)
               \end{eqnarray*} 
and         
                  \begin{eqnarray*}
                  \left( \begin{array}{c} n+M+q+a\\ 2q+a \end{array} \right)  &=&
                  \frac{ 2q+1+a}{ 2n+2M+1+a }  \left\{ 2 \left( \begin{array}{c} n+M+q+ 1+a\\2q+1+a
                  \end{array} \right)  -   \left( \begin{array}{c}
                 n+M+q+a\\2q+a \end{array} \right) \right\},
                 \end{eqnarray*}
to derive the relations:
                 \begin{subequations}
                  \begin{eqnarray}   
                  \sum_{n=-M}^N  \left. \frac{x(x+2s-3)!!}{(x-2s+1)!!}\right| _{x=2n+1+2M} 
                 &=&  \frac{1}{2(2s)}~ \left. \frac{ X(X+2s-2)!!}{(X-2s)!!} \right|_{X=2N+2+2M}; \\
                  \sum_{n=-M}^N \left.\frac{x(x+2s-2)!!}{(x-2s)!!} \right| _{x=2n+2+2M}
                 &=& \frac{1}{2(2s+1)}~ \left. \frac{ X(X+2s-1)!!}{(X-2s-1)!!} \right|_{X=2N+3+2M};\\                      
                  \sum_{n=-M+1}^N  \left. \frac{(x+2s-2)!!}{(x-2s)!!} \right| _{x=2n+2M}
                 &=&  \frac{1}{2(2s)}~ \left. \frac{ (X+2s-1)!!}{(X-2s-1)!!} \right|_{X=2N+1+2M} ;\\
                 \sum_{n=-M}^N  \left.\frac{(x+2s-1)!!}{(x-2s-1)!!}\right|_{x=2n+1+2M} 
                 &=& \frac{1}{2(2s+1)}~ \left. \frac{ (X+2s)!!}{(X-2s-2)!!} \right|_{X=2N+2+2M}.    
                 \end{eqnarray}
                 \end{subequations}
 
With the sets of relations (B4) and (B6) it is straightforward to show that the expressions 
in Theorem 2 for the hypersums satisfy condition (B2) above.  For example,  with 
$x \equiv 2n+1+2M$ and $X\equiv 2N+2+2M$, 
        \begin{eqnarray}
        \sum_{n=1}^N S^{(2M)}_{2p+1}(n) &=&  \sum_{n=-M}^N S^{(2M)}_{2p+1}(n) \nonumber \\
        &=&  \frac{(2p+1)!}{2^{2p+2+2M}} \sum_{n=-M}^N \sum_{q=0}^{p+M}
                     \left\{ C_{p+M-q}^{(2M)} ~\frac{x^{2q+2}-1}{(2q+2)!} - \frac{1}{(2q)!} \sum_{s=1}^{M} 
               \frac{C^{(2M-2s)}_{p+M+1-q-s}}{(2s)!} ~\frac{(x+2s-1)!!}{(x-2s-1)!!} \right\}\nonumber  \\
          &=&  \frac{(2p+1)!}{2^{2p+3+2M}}  \sum_{q=0}^{p+M+1}  \left\{  C^{(2M+1)}_{p+M+1-q}~ \frac{ X^{2q+1} -2^{2q}X}{(2q+1)!} 
         - \frac{1}{(2q)!} \sum_{s=1}^{M}  \frac{C^{(2M-2s)}_{p+M+1-q-s}}{(2s+1)!}   ~\frac{(X+2s)!!}{(X-2s-2)!!}  \right\} \nonumber \\
         &=&  \frac{(2p+1)!}{2^{2p+3+2M}}  \sum_{q=0}^{p+M+1}  \left\{  C^{(2M+1)}_{p+M+1-q}~ \frac{ X^{2q+1}}{(2q+1)!} 
         - \frac{1}{(2q)!} \sum_{s=0}^{M}  \frac{C^{(2M-2s)}_{p+M+1-q-s}}{(2s+1)!}   ~\frac{(X+2s)!!}{(X-2s-2)!!}  \right\} \nonumber \\
         &=&  \frac{(2p+1)!}{2^{2p+3+2M}}  \sum_{q=0}^{p+M+1}  \left\{  C^{(2M+1)}_{p+M+1-q}~ \frac{ X^{2q+1}}{(2q+1)!} 
         - \frac{1}{(2q)!} \sum_{s=1}^{M+1}  \frac{C^{(2M+2-2s)}_{p+M+2-q-s}}{(2s-1)!}   ~\frac{(X+2s-2)!!}{(X-2s)!!}  \right\} \nonumber \\
         &=& S^{(2M+1)}_{2p+1}(N),    
          \end{eqnarray}
 where we've used the relation                    
            \begin{eqnarray}
            \sum_{q=0}^p \frac{2^{2q}C^{(L+2)}_{p-q}}{(2q+1)!} &=& \sum_{q=0}^p
             \frac{C^{(L+1)}_{p-q}} {(2q)!}  
             \end{eqnarray}
from relations (\ref{y}) to move the $2^{2q}X$ term in the first summation to the second summation where 
it becomes the $s=0$ term, and then relabeled the index, $s \rightarrow s-1$.\\                 
QED

~\\
~\\
2000 Mathematics Subject Classification:  Primary:  11C08; Secondary:   11B68, 11C20, 11Y35.\\
Keywords:  Hypersums, power sums, determinant formulas, Bernoulli numbers,  Euler numbers. \\
(Concerned with sequences A036280, A036281.)

\end{document}